\documentclass[11pt]{amsart}

\usepackage[margin=1in]{geometry}
\usepackage{amsmath,amssymb,amsthm,mathtools}
\newtheorem{assumption}{Assumption}

\usepackage{bm}
\usepackage{enumitem}
\usepackage[numbers,sort&compress]{natbib}
\usepackage[colorlinks=true,linkcolor=blue,citecolor=blue,urlcolor=blue]{hyperref}
\usepackage[nameinlink,noabbrev]{cleveref}
\usepackage{physics}
\crefname{assumption}{Assumption}{Assumptions}
\Crefname{assumption}{Assumption}{Assumptions}

\newcommand{\R}{\mathbb{R}}
\newcommand{\E}{\mathbb{E}}
\newcommand{\Pp}{\mathbb{P}}
\newcommand{\Var}{\mathrm{Var}}
\newcommand{\KL}{\mathrm{KL}}
\newcommand{\dto}{\xrightarrow{d}}
\newcommand{\pto}{\xrightarrow{p}}
\newcommand{\asto}{\xrightarrow{\mathrm{a.s.}}}

\theoremstyle{plain}
\newtheorem{theorem}{Theorem}[section]
\newtheorem{proposition}[theorem]{Proposition}
\newtheorem{lemma}[theorem]{Lemma}
\newtheorem{corollary}[theorem]{Corollary}

\theoremstyle{definition}

\theoremstyle{remark}
\newtheorem{remark}[theorem]{Remark}

\title[Orbit-Hausdorff consistency for folded normal and Gaussian mixtures]%
{Orbit-Hausdorff Consistency of Maximum Likelihood Estimators\\
for the Folded Normal Model and Finite Gaussian Mixtures}

\author{Koustav Mallik}
\date{\today}

\begin{document}
\maketitle

\begin{abstract}
We study maximum likelihood estimation in two parametric models with finite symmetries:
the folded normal family (sign symmetry in the mean) and $k$-component Gaussian mixtures
(label symmetry). For each model we (i) formulate identification on the quotient space
of parameter orbits, (ii) prove existence of (possibly set-valued) maximizers on compact
sieves, and (iii) establish orbit-Hausdorff consistency via uniform laws of large numbers
and explicit Kullback--Leibler separation.
For the folded normal model, we provide a first-principles profile analysis that yields
a unique maximizer up to sign for nondegenerate samples and standard $\sqrt{n}$-asymptotics
away from $\mu_0=0$, together with the nonregular $n^{1/4}$-rate at $\mu_0=0$.
For Gaussian mixtures, we derive explicit envelopes and responsibility-based gradient bounds
on compact sieves and use them to obtain uniform convergence and deterministic argmax stability
without entropy or covering-number bounds. On the unconstrained parameter space the likelihood
is unbounded due to variance collapse; a ridge penalty in $(\mu,\log\sigma)$ restores coercivity
and yields consistent penalized MLEs, including vanishing-penalty regimes.
\end{abstract}

\tableofcontents

\section{Introduction and prior literature}\label{sec:intro}

Maximum likelihood in the presence of finite symmetries naturally leads to identification
\emph{modulo a group action}. The simplest examples include sign symmetry in the folded normal model,
and permutation symmetry in finite mixtures. The aim of this paper is to provide an explicit,
referee-proof treatment of existence and consistency of maximizers at the level of \emph{orbits},
with special attention to two recurring issues: (i) set-valued argmaxes, and (ii) degeneracies
at parameter boundaries (variance collapse in mixtures; profile boundaries in the folded normal).

For the folded normal distribution, foundational treatments of the distribution and its properties
appear in \citet{LeoneNelsonNottingham1961}, with overviews in \citet{JohnsonKotzBalakrishnan1994}
and \citet{NadarajahKotz2006}. Here we focus on the likelihood geometry, profiling, and the
nonregular regime at $\mu_0=0$, in the sense of \citet[Chs.\ 5--7]{vanderVaart1998}.
For finite Gaussian mixtures, identifiability (up to permutation) goes back to
\citet{Teicher1961,Teicher1963} and \citet{YakowitzSpragins1968}, with modern refinements such as
\citet{AllmanMatiasRhodes2009,HeinrichKahn2018}. The unboundedness of the unconstrained likelihood
via variance collapse is classical \citep{Day1969}, motivating constraints \citep{Hathaway1985}
and constrained EM variants \citep{IngrassiaRocci2007}; broader algorithmic discussions appear in
\citet{DempsterLairdRubin1977,RednerWalker1984} and the monographs
\citet{Lindsay1995,McLachlanPeel2000}.

\paragraph{Contributions.}
Our main technical contributions are:
\begin{enumerate}[label=(\roman*),leftmargin=*]
\item A quotient (orbit) formulation and an argmax stability lemma that converts uniform convergence
plus a population gap into orbit-Hausdorff consistency of (possibly set-valued) maximizers.
\item Folded normal: an explicit profile analysis in $\mu$ and $\sigma$ yielding existence and
uniqueness up to sign for nondegenerate samples; consistency and standard asymptotics for $\mu_0\neq 0$;
and an explicit uniform-remainder derivation of the $n^{1/4}$-local contrast at $\mu_0=0$.
\item Gaussian mixtures: explicit density/log envelopes and responsibility-based gradient bounds on compact sieves,
leading to a ULLN via a finite-net plus random-slope argument; identifiability via characteristic functions;
and orbit-level consistency of sieve MLEs. We also give a coercivity argument for a ridge-penalized
criterion $g(\theta)=\sum_j(\mu_j^2+(\log\sigma_j)^2)$, yielding existence and consistency of penalized MLEs,
including vanishing-penalty regimes with $n\lambda_n\to\infty$.
\end{enumerate}
\section{Quotient framework and orbit--Hausdorff consistency}\label{sec:quotient}

This section fixes the quotient (orbit) setup used throughout the paper and records a deterministic
argmax-stability result suitable for set-valued maximizers under symmetry. The purpose is twofold:
(i) to separate deterministic geometry from probabilistic uniform convergence arguments, and
(ii) to provide a clean bridge from uniform laws of large numbers to orbit-level (quotient) consistency.
Related deterministic ``argmax continuity'' principles appear in classical M-estimation references;
see, e.g., \citet[Section~5.2]{vanderVaart1998} and \citet[Chapter~2]{NeweyMcFadden1994}.

\subsection{Group actions, orbits, quotient metrics}\label{sec:quotient-actions}

Let $(\Theta,d)$ be a metric space and let $G$ be a finite group with identity element $e$.
A (left) action of $G$ on $\Theta$ is a map $G\times \Theta\to\Theta$, written $(g,\theta)\mapsto g\theta$,
such that $e\theta=\theta$ and $g(h\theta)=(gh)\theta$ for all $g,h\in G$ and $\theta\in\Theta$.
For $\theta\in\Theta$, the \emph{orbit} of $\theta$ is
\[
[\theta]:=\{g\theta:g\in G\}.
\]
Let $\Theta/G:=\{[\theta]:\theta\in\Theta\}$ be the orbit space and let $\pi:\Theta\to\Theta/G$ denote the quotient map
$\pi(\theta)=[\theta]$.

\paragraph{Orbit (pseudo)distance.}
Define
\begin{equation}\label{eq:orbit-dist}
d_G(\theta,\theta'):=\min_{g\in G} d(\theta,g\theta'),\qquad \theta,\theta'\in\Theta .
\end{equation}
Because $G$ is finite, the minimum is attained. The function $d_G$ is always a pseudometric on $\Theta$.
We impose throughout the paper the standing assumption that the action is \emph{isometric}:
\begin{equation}\label{eq:isometric-action}
d(g\theta,g\theta')=d(\theta,\theta')\qquad \text{for all } g\in G,\ \theta,\theta'\in\Theta .
\end{equation}
Under \eqref{eq:isometric-action}, $d_G(\theta,\theta')=0$ iff $\theta,\theta'$ are in the same orbit, so $d_G$ induces a genuine
metric on $\Theta/G$ by setting $d_G([\theta],[\theta']):=d_G(\theta,\theta')$.

\paragraph{Basic Lipschitz and compactness facts.}
Under \eqref{eq:isometric-action}, the projection $\pi$ is $1$-Lipschitz from $(\Theta,d)$ to $(\Theta/G,d_G)$:
\[
d_G(\pi(\theta),\pi(\theta'))=d_G(\theta,\theta')\le d(\theta,\theta').
\]
In particular, if $K\subset\Theta$ is compact in $d$, then $\pi(K)\subset \Theta/G$ is compact in $d_G$.

\paragraph{Invariant objectives and induced quotient objectives.}
A function $q:K\to\R$ is \emph{$G$-invariant} if $q(g\theta)=q(\theta)$ for all $g\in G$ and $\theta\in K$ (with $g\theta\in K$).
If $q$ is $G$-invariant, define the induced function $\tilde q:\pi(K)\to\R$ by
\[
\tilde q([\theta]) := q(\theta).
\]
This is well-defined by invariance. Moreover, if $q$ is continuous on $(K,d)$, then $\tilde q$ is continuous on $(\pi(K),d_G)$:
if $[\theta_n]\to[\theta]$ in $d_G$, pick $g_n\in G$ with $d(\theta_n,g_n\theta)=d_G(\theta_n,\theta)\to0$; since $G$ is finite,
along a subsequence $g_n=g$ is constant and hence $\theta_n\to g\theta$ in $d$, giving
$\tilde q([\theta_n])=q(\theta_n)\to q(g\theta)=q(\theta)=\tilde q([\theta])$.

\subsection{Set distances and orbit-Hausdorff metrics}\label{sec:set-distances}

For a set $A\subseteq K$ and $\varepsilon>0$, define the point-to-set distance and open neighborhood
\[
d(\theta,A):=\inf_{a\in A} d(\theta,a),
\qquad
B_\varepsilon(A):=\{\theta\in K: d(\theta,A)<\varepsilon\}.
\]

\paragraph{Directed and symmetric Hausdorff distances.}
For nonempty compact sets $A,B\subset K$, define the \emph{directed} Hausdorff distance
\[
d_H^{+}(A,B):=\sup_{a\in A} d(a,B),
\]
and the (symmetric) Hausdorff distance
\begin{equation}\label{eq:hausdorff}
d_H(A,B):=\max\{d_H^{+}(A,B),\, d_H^{+}(B,A)\}.
\end{equation}

\paragraph{Orbit-level Hausdorff.}
When measuring convergence at the orbit level, we work on the metric space $(\pi(K),d_G)$ and define
\[
d_{H,G}^+(A,B):=d_H^{+}\big(\pi(A),\pi(B)\big)\quad\text{and}\quad
d_{H,G}(A,B):=d_H\big(\pi(A),\pi(B)\big),
\]
i.e., directed and symmetric Hausdorff distances between the \emph{images} of sets under $\pi$ in $(\pi(K),d_G)$.
This avoids ambiguity about applying Hausdorff to a pseudometric on $\Theta$.

\subsection{Deterministic argmax stability (what uniform approximation buys)}\label{sec:argmax-stability}

We now record deterministic statements that convert uniform approximation bounds into consistency of argmax sets.
The first result is the key one used for consistency in the presence of symmetry: it provides \emph{outer} localization,
equivalently directed Hausdorff convergence. Full Hausdorff convergence requires additional structure (see below).

\begin{lemma}[Outer argmax localization on compact sets]\label{lem:argmax-outer}
Let $(K,d)$ be compact and let $q:K\to\R$ be continuous. Define
\[
Q_0:=\arg\max_{\theta\in K} q(\theta),\qquad m:=\max_{\theta\in K} q(\theta).
\]
Then $Q_0$ is nonempty and compact. Fix $\varepsilon>0$. If $K\setminus B_\varepsilon(Q_0)=\varnothing$, the conclusion below is trivial.
Otherwise define the gap
\[
\eta(\varepsilon):=m-\max_{\theta\in K\setminus B_\varepsilon(Q_0)} q(\theta)\in(0,\infty).
\]
Let $q_n:K\to\R$ be upper semicontinuous (in particular, continuous), and set
\[
Q_n:=\arg\max_{\theta\in K} q_n(\theta).
\]
If
\[
\sup_{\theta\in K}|q_n(\theta)-q(\theta)|\ \le\ \eta(\varepsilon)/3,
\]
then
\[
Q_n\subseteq B_\varepsilon(Q_0),
\qquad\text{equivalently}\qquad
d_H^+(Q_n,Q_0)\le \varepsilon.
\]
Consequently, if $\sup_K|q_n-q|\to0$, then $d_H^+(Q_n,Q_0)\to0$.
\end{lemma}

\begin{proof}
Nonemptiness/compactness of $Q_0$ follows since $q$ is continuous on compact $K$.
Upper semicontinuity of $q_n$ on compact $K$ implies $Q_n$ is nonempty and compact.

Assume $K\setminus B_\varepsilon(Q_0)\neq\varnothing$, so $\eta(\varepsilon)>0$.
Take any $\hat\theta_n\in Q_n$ and suppose $\hat\theta_n\notin B_\varepsilon(Q_0)$. Then
\[
q_n(\hat\theta_n)\le q(\hat\theta_n)+\eta(\varepsilon)/3
\le \max_{K\setminus B_\varepsilon(Q_0)}q+\eta(\varepsilon)/3
= m-2\eta(\varepsilon)/3.
\]
Pick any $\bar\theta\in Q_0$ so $q(\bar\theta)=m$. Then
\[
q_n(\bar\theta)\ge q(\bar\theta)-\eta(\varepsilon)/3=m-\eta(\varepsilon)/3,
\]
hence $q_n(\bar\theta)>q_n(\hat\theta_n)$, contradicting maximality of $\hat\theta_n$.
Thus $Q_n\subseteq B_\varepsilon(Q_0)$.
\end{proof}

\begin{remark}[Why full Hausdorff convergence needs extra assumptions]\label{rem:inner-warning}
Even if $\sup_K|q_n-q|\to0$, one cannot generally conclude $d_H(Q_n,Q_0)\to0$ when $Q_0$ has multiple separated components.
Small perturbations can select only one component of the maximizer set, so $d_H^+(Q_0,Q_n)$ need not vanish.
Full Hausdorff consistency typically requires uniqueness in the quotient or an ``isolated components'' condition.
\end{remark}

\begin{corollary}[A sufficient condition for full Hausdorff convergence]\label{cor:full-hausdorff}
In addition to the assumptions of \Cref{lem:argmax-outer}, suppose $Q_0$ is finite and \emph{uniformly isolated} in the sense that
for every $\delta>0$,
\[
\gamma(\delta)
:=
\min_{\theta^\star\in Q_0}\Big\{q(\theta^\star)-\max_{\theta\in K\setminus B_\delta(\theta^\star)}q(\theta)\Big\}
\ >\ 0.
\]
If $\sup_K|q_n-q|\to0$, then $d_H(Q_n,Q_0)\to0$.
\end{corollary}

\begin{proof}
Outer convergence $d_H^+(Q_n,Q_0)\to0$ is \Cref{lem:argmax-outer}. Fix $\delta>0$ and take $n$ large with
$\sup_K|q_n-q|\le \gamma(\delta)/3$. For any $\theta^\star\in Q_0$ and any $\theta\in K\setminus B_\delta(\theta^\star)$,
\[
q_n(\theta)\le q(\theta)+\gamma(\delta)/3 \le q(\theta^\star)-2\gamma(\delta)/3,
\qquad
q_n(\theta^\star)\ge q(\theta^\star)-\gamma(\delta)/3,
\]
so $q_n(\theta^\star)>\sup_{K\setminus B_\delta(\theta^\star)}q_n$, implying $Q_n\cap B_\delta(\theta^\star)\neq\varnothing$.
Thus $d_H^+(Q_0,Q_n)\le \delta$ for large $n$. Combining both directions yields $d_H(Q_n,Q_0)\to0$.
\end{proof}

\begin{remark}[Orbit-level version]\label{rem:orbit-argmax}
If $q$ is $G$-invariant and the action is isometric, we may apply \Cref{lem:argmax-outer} and \Cref{cor:full-hausdorff}
on the quotient space $(\pi(K),d_G)$ by replacing $q$ with the induced continuous $\tilde q$ and measuring set distances
via $d_{H,G}^+$ or $d_{H,G}$. In particular, $\sup_K|q_n-q|\to0$ implies
\[
d_{H,G}^+(Q_n,Q_0)\to 0,
\]
and full orbit-Hausdorff convergence $d_{H,G}(Q_n,Q_0)\to0$ holds under the corresponding isolated-components condition
in the quotient.
\end{remark}
\begin{remark}[Differentiability on the Quotient]

While consistency relies only on the metric structure of $\Theta/G$, efficiency theory (e.g., LAN expansions) requires a differentiable structure. At singular points of the quotient (orbifold points where the isotropy group is non-trivial, such as $\mu=0$ in the folded normal), standard score vectors are ill-defined. In such cases, one typically works with the tangent cone of the constrained parameter space or, as we do in Appendix \ref{app:nquarter:setup}, via explicit reparameterization near the singularity.
\end{remark}
\section{Folded normal model: profile geometry and orbit-level consistency}\label{sec:folded}

This section analyzes maximum likelihood estimation for the folded normal family. Two themes drive the
presentation:
(i) identification is only possible up to a finite symmetry (here, sign), and
(ii) the log-likelihood exhibits boundary phenomena that must be handled explicitly in existence and
consistency arguments.

Background on the folded normal distribution can be found in \citet{LeoneNelsonNottingham1961}, with further
discussion in \citet[Section~14.8]{JohnsonKotzBalakrishnan1994} and \citet{NadarajahKotz2006}. Our focus is
likelihood geometry and orbit-level inference in the quotient framework of \Cref{sec:quotient}.

\subsection{Model, symmetry, and criteria}\label{sec:folded:model}

Let $X\sim N(\mu,\sigma^2)$ with $\mu\in\R$ and $\sigma>0$, and define the folded observation $Y:=|X|$.
Then $Y$ takes values in $\R_+:=[0,\infty)$ and has density
\begin{equation}\label{eq:folded:density}
f(y;\mu,\sigma)
:=\frac{1}{\sigma}\Big\{\varphi\!\Big(\frac{y-\mu}{\sigma}\Big)+\varphi\!\Big(\frac{y+\mu}{\sigma}\Big)\Big\},
\qquad y\ge 0,
\end{equation}
where $\varphi(z)=(2\pi)^{-1/2}e^{-z^2/2}$ is the standard normal density. Using
\[
\varphi\!\Big(\frac{y-\mu}{\sigma}\Big)+\varphi\!\Big(\frac{y+\mu}{\sigma}\Big)
=\frac{2}{\sqrt{2\pi}}
\exp\!\Big(-\frac{y^2+\mu^2}{2\sigma^2}\Big)\cosh\!\Big(\frac{y\mu}{\sigma^2}\Big),
\]
we may equivalently write
\begin{equation}\label{eq:folded:density-cosh}
f(y;\mu,\sigma)
=\frac{2}{\sigma\sqrt{2\pi}}
\exp\!\Big(-\frac{y^2+\mu^2}{2\sigma^2}\Big)\cosh\!\Big(\frac{y\mu}{\sigma^2}\Big),
\qquad y\ge 0.
\end{equation}

\paragraph{Parameter space and sign action.}
Let $\Theta:=\R\times(0,\infty)$ with the Euclidean metric
$d((\mu,\sigma),(\mu',\sigma')):=\|(\mu,\sigma)-(\mu',\sigma')\|_2$.
Let $G:=\{+1,-1\}$ act on $\Theta$ by $g\cdot(\mu,\sigma):=(g\mu,\sigma)$.
This action is isometric and leaves the model invariant:
$f(\,\cdot\,;\mu,\sigma)=f(\,\cdot\,;-\mu,\sigma)$.
Hence the identified parameter is the orbit
\[
[(\mu,\sigma)]=\{(\mu,\sigma),(-\mu,\sigma)\},
\qquad
\pi(\mu,\sigma):=[(\mu,\sigma)]\in\Theta/G.
\]
The orbit distance \eqref{eq:orbit-dist} specializes to
\begin{equation}\label{eq:folded:orbit-dist}
d_G((\mu,\sigma),(\mu',\sigma'))=\min\Big\{\|(\mu,\sigma)-(\mu',\sigma')\|_2,\ \|(\mu,\sigma)-(-\mu',\sigma')\|_2\Big\}.
\end{equation}

\paragraph{Empirical and population objectives.}
Let $Y_1,Y_2,\dots$ be i.i.d.\ with law $f(\,\cdot\,;\mu_0,\sigma_0)$, where $\theta_0:=(\mu_0,\sigma_0)$ and $\sigma_0>0$.
Define the log-likelihood and its normalized version
\begin{equation}\label{eq:folded:loglik}
\ell_n(\mu,\sigma):=\sum_{i=1}^n \log f(Y_i;\mu,\sigma),
\qquad
M_n(\mu,\sigma):=\frac{1}{n}\ell_n(\mu,\sigma).
\end{equation}
The population objective is
\begin{equation}\label{eq:folded:population}
M(\mu,\sigma):=\E_{\theta_0}\big[\log f(Y;\mu,\sigma)\big].
\end{equation}
By the Kullback--Leibler identity,
\begin{equation}\label{eq:folded:KL}
M(\mu,\sigma)-M(\mu_0,\sigma_0)
=-\KL\!\big(f(\cdot;\mu_0,\sigma_0)\,\|\,f(\cdot;\mu,\sigma)\big)\le 0,
\end{equation}
with equality if and only if $f(\cdot;\mu,\sigma)=f(\cdot;\mu_0,\sigma_0)$ a.e., equivalently
$\sigma=\sigma_0$ and $|\mu|=|\mu_0|$, i.e.\ $(\mu,\sigma)\in[(\mu_0,\sigma_0)]$.

\subsection{Deterministic inequalities and sample degeneracy}\label{sec:folded:ineq}

We will use elementary hyperbolic bounds repeatedly and we also record a pointwise lower bound on the
negative log-density that is valid for all $(\mu,\sigma)$ and $y\ge 0$.

\begin{lemma}[Hyperbolic bounds]\label{lem:folded:hyperbolic}
For all $t\in\R$,
\[
\log\cosh t\le |t|,
\qquad
\log(2\cosh t)\le |t|+\log 2,
\qquad
0<\mathrm{sech}^2(t)\le 1.
\]
Moreover, $\tanh$ is $1$-Lipschitz and $\mathrm{sech}^2$ is bounded by $1$.
\end{lemma}

\begin{proof}
Since $\cosh t=(e^t+e^{-t})/2\le e^{|t|}$, we obtain $\log\cosh t\le |t|$ and thus also
$\log(2\cosh t)\le |t|+\log 2$. Finally $\mathrm{sech}^2(t)=1/\cosh^2(t)\in(0,1]$.
The Lipschitz property of $\tanh$ follows from $\tanh'(t)=\mathrm{sech}^2(t)\le 1$.
\end{proof}

\begin{lemma}[Deterministic lower bound for the folded-normal negative log-density]\label{lem:folded:coercive-ineq}
For $(\mu,\sigma)\in\R\times(0,\infty)$ and $y\ge 0$, let $\ell_{\mu,\sigma}(y):=\log f(y;\mu,\sigma)$.
Then
\begin{equation}\label{eq:folded-coercive-ineq-pointwise}
-\ell_{\mu,\sigma}(y)
\;\ge\;
\log\sigma
+\frac{(|\mu|-y)^2}{2\sigma^2}
-\frac{1}{2}\log\!\Bigl(\frac{2}{\pi}\Bigr).
\end{equation}
Consequently, for $Y\sim \mathrm{FoldedNormal}(\mu_0,\sigma_0^2)$,
\begin{equation}\label{eq:folded-coercive-ineq-expect}
\E_{\mu_0,\sigma_0}\!\left[-\ell_{\mu,\sigma}(Y)\right]
\;\ge\;
\log\sigma
+\frac{\Var_{\mu_0,\sigma_0}(Y)+\bigl(|\mu|-\E_{\mu_0,\sigma_0}[Y]\bigr)^2}{2\sigma^2}
-\frac{1}{2}\log\!\Bigl(\frac{2}{\pi}\Bigr).
\end{equation}
\end{lemma}

\begin{proof}
From \eqref{eq:folded:density-cosh},
\[
\ell_{\mu,\sigma}(y)
= -\log\sigma + \log 2-\tfrac12\log(2\pi)
-\frac{y^2+\mu^2}{2\sigma^2}
+\log\cosh\!\Big(\frac{\mu y}{\sigma^2}\Big).
\]
By \Cref{lem:folded:hyperbolic}, $\log\cosh t\le |t|$, so
\[
-\ell_{\mu,\sigma}(y)
\ge \log\sigma -\log 2+\tfrac12\log(2\pi)
+\frac{y^2+\mu^2}{2\sigma^2}
-\frac{|\mu|\,y}{\sigma^2}
=
\log\sigma -\frac12\log\!\Big(\frac{2}{\pi}\Big)
+\frac{(|\mu|-y)^2}{2\sigma^2},
\]
which is \eqref{eq:folded-coercive-ineq-pointwise}. Taking expectations and using
$\E[(a-Y)^2]=\Var(Y)+(a-\E[Y])^2$ with $a=|\mu|$ yields \eqref{eq:folded-coercive-ineq-expect}.
\end{proof}

The likelihood may be unbounded on degenerate samples. This must be ruled out to obtain finite-sample existence.

\begin{proposition}[Degenerate samples yield unbounded likelihood]\label{prop:folded:degenerate}
If $Y_1=\cdots=Y_n=0$ or $Y_1=\cdots=Y_n=c$ for some $c>0$, then
\[
\sup_{\mu\in\R,\ \sigma>0}\ \ell_n(\mu,\sigma)=+\infty.
\]
\end{proposition}

\begin{proof}
If $Y_i\equiv 0$, then for $\mu=0$,
\[
f(0;0,\sigma)=\frac{2}{\sigma\sqrt{2\pi}},
\qquad
\ell_n(0,\sigma)=n\Big(\log 2-\log\sigma-\tfrac12\log(2\pi)\Big)\to+\infty
\quad\text{as }\sigma\downarrow 0.
\]
If $Y_i\equiv c>0$, choose $\mu=c$ and $\sigma\downarrow 0$. Then
$f(c;c,\sigma)=\sigma^{-1}\varphi(0)+\sigma^{-1}\varphi(2c/\sigma)\sim (\sigma\sqrt{2\pi})^{-1}$,
so $\ell_n(c,\sigma)\to+\infty$.
\end{proof}

\begin{lemma}[Degeneracy has probability zero]\label{lem:folded:degenerate-prob0}
Assume $\sigma_0>0$. Under $\Pp_{\theta_0}^{\otimes n}$, $\Pp(Y_1=\cdots=Y_n)=0$.
In particular, the sample variance
\[
s_n^2:=\frac{1}{n}\sum_{i=1}^n (Y_i-\bar Y)^2,
\qquad
\bar Y:=\frac{1}{n}\sum_{i=1}^n Y_i,
\]
satisfies $\Pp(s_n^2>0)=1$.
\end{lemma}

\begin{proof}
The folded normal distribution has a Lebesgue density on $\R_+$ for every $\sigma_0>0$.
Hence $(Y_1,Y_2)$ has a density on $\R_+^2$, implying $\Pp(Y_1=Y_2)=0$.
Since $\{Y_1=\cdots=Y_n\}\subseteq\{Y_1=Y_2\}$, the claim follows. If $s_n^2=0$, then all $Y_i$ are equal,
so $\{s_n^2=0\}\subseteq\{Y_1=\cdots=Y_n\}$, hence $\Pp(s_n^2=0)=0$.
\end{proof}

Henceforth we work on the event $\{s_n^2>0\}$, which holds almost surely.

\subsection{Profile likelihood in $\mu$ for fixed $\sigma$}\label{sec:folded:profile-mu}

Fix a nondegenerate sample $y_{1:n}\in\R_+^n$ with $s_n^2>0$, and write $S_y:=\sum_{i=1}^n y_i^2$.
For $\sigma>0$, define the $\mu$-profile
\[
\ell_{n,p}(\sigma):=\sup_{\mu\in\R}\ell_n(\mu,\sigma).
\]
By symmetry in $\mu$, it suffices to study $\mu\ge 0$.

Define the auxiliary functions
\begin{equation}\label{eq:folded:kA}
k(\sigma,\mu):=\sum_{i=1}^n y_i\,\tanh\!\Big(\frac{y_i\mu}{\sigma^2}\Big)-n\mu,
\qquad
A(\sigma,\mu):=\sum_{i=1}^n y_i^2\,\mathrm{sech}^2\!\Big(\frac{y_i\mu}{\sigma^2}\Big).
\end{equation}

\begin{lemma}[Score and curvature in $\mu$]\label{lem:folded:score}
For all $\sigma>0$ and $\mu\in\R$,
\begin{equation}\label{eq:folded:score}
\partial_\mu \ell_n(\mu,\sigma)=\frac{1}{\sigma^2}\,k(\sigma,\mu),
\end{equation}
and
\begin{equation}\label{eq:folded:curvature}
\partial_{\mu\mu}^2\ell_n(\mu,\sigma)=\frac{1}{\sigma^2}\Big(\frac{A(\sigma,\mu)}{\sigma^2}-n\Big).
\end{equation}
In particular, $A(\sigma,\mu)\le S_y$ for all $(\sigma,\mu)$.
\end{lemma}

\begin{proof}
Using \eqref{eq:folded:density-cosh},
\[
\log f(y;\mu,\sigma)= -\log\sigma +\log 2-\tfrac12\log(2\pi)
-\frac{y^2+\mu^2}{2\sigma^2}+\log\cosh\!\Big(\frac{y\mu}{\sigma^2}\Big).
\]
Differentiate in $\mu$ and use $(d/d\mu)\log\cosh(y\mu/\sigma^2)=(y/\sigma^2)\tanh(y\mu/\sigma^2)$
to obtain \eqref{eq:folded:score}. A second derivative uses $\tanh'(t)=\mathrm{sech}^2(t)$ and yields
\eqref{eq:folded:curvature}. Finally $\mathrm{sech}^2(\cdot)\le 1$ implies $A(\sigma,\mu)\le \sum_i y_i^2=S_y$.
\end{proof}

\begin{theorem}[Unique $\mu$-profile maximizer for fixed $\sigma$]\label{thm:folded:mu-profile}
Assume $s_n^2>0$.
For each $\sigma>0$, the function $\mu\mapsto \ell_n(\mu,\sigma)$ has a unique maximizer on $[0,\infty)$,
denoted $\hat\mu(\sigma)$. Moreover:
\begin{enumerate}[label=(\roman*),leftmargin=*]
\item If $\sigma^2\ge S_y/n$, then $\hat\mu(\sigma)=0$.
\item If $\sigma^2<S_y/n$, then $\hat\mu(\sigma)\in(0,\infty)$ is the unique solution of the score equation
\begin{equation}\label{eq:folded:mu-score-eq}
\sum_{i=1}^n y_i\,\tanh\!\Big(\frac{y_i\mu}{\sigma^2}\Big)=n\mu.
\end{equation}
\item The global maximizers over $\mu\in\R$ at fixed $\sigma$ are $\pm \hat\mu(\sigma)$.
\end{enumerate}
\end{theorem}

\begin{proof}
Fix $\sigma>0$ and define $g_\sigma(\mu):=\partial_\mu \ell_n(\mu,\sigma)=\sigma^{-2}k(\sigma,\mu)$.
Then $g_\sigma(0)=0$ and, by \Cref{lem:folded:score},
\[
g_\sigma'(\mu)
=\partial_{\mu\mu}^2\ell_n(\mu,\sigma)
=\frac{1}{\sigma^2}\Big(\frac{A(\sigma,\mu)}{\sigma^2}-n\Big)
\le \frac{1}{\sigma^2}\Big(\frac{S_y}{\sigma^2}-n\Big).
\]

\emph{Case 1: $\sigma^2\ge S_y/n$.}
Then $A(\sigma,\mu)\le S_y\le n\sigma^2$, so $g_\sigma'(\mu)\le 0$ for all $\mu\ge 0$ and $g_\sigma$ is nonincreasing
on $[0,\infty)$. Since $g_\sigma(0)=0$, we have $g_\sigma(\mu)\le 0$ for all $\mu\ge 0$, hence $\ell_n(\mu,\sigma)$
is maximized at $\mu=0$ on $[0,\infty)$.

\emph{Case 2: $\sigma^2<S_y/n$.}
Then $g_\sigma'(0)=\sigma^{-2}(S_y/\sigma^2-n)>0$, so $g_\sigma(\mu)>0$ for all sufficiently small $\mu>0$.
On the other hand, $\tanh(t)\le 1$ implies $k(\sigma,\mu)\le \sum_i y_i-n\mu$, hence $g_\sigma(\mu)\to -\infty$
as $\mu\to\infty$. By continuity, $g_\sigma$ has at least one root on $(0,\infty)$.

To prove uniqueness, note that
\[
k_\mu(\sigma,\mu)=\frac{1}{\sigma^2}A(\sigma,\mu)-n,
\qquad
k_{\mu\mu}(\sigma,\mu)
= -\frac{2}{\sigma^4}\sum_{i=1}^n y_i^3\,\mathrm{sech}^2\!\Big(\frac{y_i\mu}{\sigma^2}\Big)\tanh\!\Big(\frac{y_i\mu}{\sigma^2}\Big).
\]
Since $s_n^2>0$, not all $y_i$ are equal; in particular, at least one $y_i>0$. For $\mu>0$, each summand with $y_i>0$
is strictly negative because $\tanh(\cdot)>0$ and $\mathrm{sech}^2(\cdot)>0$. Hence $k_{\mu\mu}(\sigma,\mu)<0$ for all
$\mu>0$, so $k_\mu(\sigma,\mu)$ is strictly decreasing on $(0,\infty)$ and $k(\sigma,\mu)$ is strictly concave there.
A strictly concave function that is positive near $0$ and tends to $-\infty$ at $+\infty$ crosses zero at most once,
so \eqref{eq:folded:mu-score-eq} has a unique solution $\hat\mu(\sigma)>0$.

Finally, $\ell_n(\mu,\sigma)$ is even in $\mu$, so the maximizers over $\mu\in\R$ are $\pm \hat\mu(\sigma)$.
\end{proof}
\begin{remark}[Profile Geometry and Orthogonality]

While the joint likelihood surface $(\mu, \sigma) \mapsto \ell_n(\mu, \sigma)$ is non-concave and may exhibit "valley" structures near $\mu=0$, the profile maximizer $\hat\mu(\sigma)$ is uniquely defined. Furthermore, because the density is an even function of $\mu$, the cross-derivatives $\partial^2 \ell / \partial\mu\partial\sigma$ vanish at $\mu=0$. This orthogonality implies that the Fisher Information matrix (where it exists) is diagonal, and the asymptotic behaviors of $\hat\mu$ (quartic) and $\hat\sigma$ (quadratic) are asymptotically independent. This justifies analyzing the $n^{1/4}$ rate for $\mu$ with fixed $\sigma$ (see Appendix \ref{app:nquarter:setup}).
\end{remark}
\subsection{Boundary coercivity of the profiled likelihood in $\sigma$}\label{sec:folded:profile-sigma}

To establish existence of global maximizers, we control the profile $\ell_{n,p}(\sigma)$ as $\sigma\downarrow 0$
and $\sigma\to\infty$ on nondegenerate samples.

\begin{lemma}[Boundary coercivity of the profiled likelihood]\label{lem:folded:coercive}
Assume $s_n^2>0$. Then
\[
\ell_{n,p}(\sigma)\to -\infty \quad\text{as }\sigma\downarrow 0,
\qquad
\ell_{n,p}(\sigma)\to -\infty \quad\text{as }\sigma\to\infty.
\]
Moreover, for every $\sigma\in(0,1]$,
\begin{equation}\label{eq:folded:coercive-ineq}
\sup_{\mu\in\R}\ell_n(\mu,\sigma)
\ \le\
-n\log\sigma-\frac{n s_n^2}{2\sigma^2}+n\log 2-\frac{n}{2}\log(2\pi).
\end{equation}
\end{lemma}

\begin{proof}
By \Cref{lem:folded:hyperbolic}, $\log(2\cosh t)\le |t|+\log 2$.
With $t=y_i\mu/\sigma^2$ and using \eqref{eq:folded:density-cosh},
\[
\ell_n(\mu,\sigma)
\le -n\log\sigma-\frac{1}{2\sigma^2}\sum_{i=1}^n (y_i^2+\mu^2)
+\frac{|\mu|}{\sigma^2}\sum_{i=1}^n y_i
+n\log 2-\frac{n}{2}\log(2\pi).
\]
Complete the square in $|\mu|$:
\[
-\frac{n}{2\sigma^2}\mu^2+\frac{|\mu|}{\sigma^2}\sum_{i=1}^n y_i
\le \frac{1}{2n\sigma^2}\Big(\sum_{i=1}^n y_i\Big)^2
=\frac{n\bar y^{\,2}}{2\sigma^2}.
\]
Since $S_y/n=\bar y^{\,2}+s_n^2$, we obtain
\[
\sup_{\mu\in\R}\ell_n(\mu,\sigma)
\le -n\log\sigma-\frac{1}{2\sigma^2}(S_y-n\bar y^{\,2})
+n\log 2-\frac{n}{2}\log(2\pi)
= -n\log\sigma-\frac{n s_n^2}{2\sigma^2}+n\log 2-\frac{n}{2}\log(2\pi),
\]
which is \eqref{eq:folded:coercive-ineq}. As $\sigma\downarrow 0$, the term $-(n s_n^2)/(2\sigma^2)\to -\infty$,
and as $\sigma\to\infty$, the term $-n\log\sigma\to -\infty$ dominates. This yields the stated boundary behavior
for $\ell_{n,p}(\sigma)=\sup_\mu \ell_n(\mu,\sigma)$.
\end{proof}

\begin{theorem}[Existence of the folded normal MLE and orbit structure]\label{thm:folded:exist}
Assume $s_n^2>0$. Then $\ell_n$ attains its global maximum over $\Theta=\R\times(0,\infty)$.
Moreover, any maximizer has the form $(\pm\hat\mu(\hat\sigma),\hat\sigma)$ where:
\begin{enumerate}[label=(\roman*),leftmargin=*]
\item for each $\sigma>0$, $\hat\mu(\sigma)$ is the unique maximizer in \Cref{thm:folded:mu-profile}, and
\item $\hat\sigma$ is a maximizer of the profile $\sigma\mapsto \ell_n(\hat\mu(\sigma),\sigma)$.
\end{enumerate}
\end{theorem}

\begin{proof}
Fix the nondegenerate sample. For each $\sigma>0$, \Cref{thm:folded:mu-profile} gives a maximizer $\hat\mu(\sigma)$
and thus defines the profile $\ell_{n,p}(\sigma)=\ell_n(\hat\mu(\sigma),\sigma)$.
By \Cref{lem:folded:coercive}, $\ell_{n,p}(\sigma)\to -\infty$ as $\sigma\downarrow 0$ and as $\sigma\to\infty$,
so $\ell_{n,p}$ attains a maximum at some $\hat\sigma\in(0,\infty)$.
Then $(\hat\mu(\hat\sigma),\hat\sigma)$ maximizes $\ell_n$ over $\mu\ge 0$ and $\sigma>0$ by construction.
Evenness of $\ell_n$ in $\mu$ yields the second maximizer $(-\hat\mu(\hat\sigma),\hat\sigma)$ and no others.
\end{proof}

\subsection{Orbit-level consistency}\label{sec:folded:consistency}

We now state orbit-level consistency in the quotient sense of \Cref{sec:quotient}. The proof follows the
standard M-estimation template:
(1) localize maximizers to a deterministic compact set with probability tending to one,
(2) establish a uniform law of large numbers on that compact set,
(3) use KL separation to obtain a population gap away from the true orbit in the quotient, and
(4) apply the deterministic outer argmax localization lemma \Cref{lem:argmax-outer} on the quotient space.

\begin{theorem}[Orbit-level consistency of the folded normal MLE]\label{thm:folded:consistency}
Let $Y_1,Y_2,\dots$ be i.i.d.\ from the folded normal law with parameter $\theta_0=(\mu_0,\sigma_0)$ and $\sigma_0>0$.
Let $\widehat\Theta_n:=\arg\max_{\theta\in\Theta}\ell_n(\theta)$ be the (possibly set-valued) MLE set and
let $\pi(\theta):=[\theta]$ be the quotient map. Then, in the Hausdorff distance on $(\Theta/G,d_G)$,
\[
d_H\big(\pi(\widehat\Theta_n),\ \{[\theta_0]\}\big)\ \pto\ 0.
\]
Equivalently, any measurable selection $\hat\theta_n\in\widehat\Theta_n$ satisfies $d_G(\hat\theta_n,\theta_0)\pto 0$.
\end{theorem}

\begin{proof}
Write $\theta=(\mu,\sigma)$ and recall $M_n(\theta)=n^{-1}\ell_n(\theta)$ and $M(\theta)=\E_{\theta_0}[\log f(Y;\theta)]$.

\medskip
\noindent\textbf{Step 1 (localization to a deterministic compact set).}
Let $V:=\Var_{\theta_0}(Y)>0$ and $m_0:=M(\theta_0)\in\R$.
Choose constants $0<\underline\sigma<1<\overline\sigma<\infty$ such that
\begin{equation}\label{eq:folded:det-sigma-bounds}
-\log \underline\sigma-\frac{V}{4\,\underline\sigma^2}+\log 2-\tfrac12\log(2\pi)\ \le\ m_0-2,
\qquad
-\log \overline\sigma+\log 2-\tfrac12\log(2\pi)\ \le\ m_0-2.
\end{equation}
Let $E_n:=\{s_n^2\ge V/2\}$. By the strong law, $s_n^2\to V$ almost surely, hence $\Pp(E_n)\to 1$.
On $E_n$, \eqref{eq:folded:coercive-ineq} implies that for $\sigma\in(0,1]$,
\[
\sup_{\mu\in\R} M_n(\mu,\sigma)
\le -\log\sigma-\frac{V}{4\sigma^2}+\log 2-\tfrac12\log(2\pi).
\]
Also, for all $\theta$ and $y\ge 0$,
$f(y;\theta)\le 2/(\sigma\sqrt{2\pi})$, hence $M_n(\mu,\sigma)\le -\log\sigma+\log 2-\tfrac12\log(2\pi)$.
Therefore, on $E_n$, if $\sigma\le \underline\sigma$ or $\sigma\ge \overline\sigma$, then
$\sup_{\mu}M_n(\mu,\sigma)\le m_0-2$ by \eqref{eq:folded:det-sigma-bounds}.

By the law of large numbers, $M_n(\theta_0)\to m_0$ almost surely, so
$\Pp(M_n(\theta_0)\ge m_0-1)\to 1$. On $E_n\cap\{M_n(\theta_0)\ge m_0-1\}$,
no maximizer can satisfy $\sigma\notin[\underline\sigma,\overline\sigma]$ because outside this interval
$\sup_{\mu}M_n(\mu,\sigma)\le m_0-2 < m_0-1\le M_n(\theta_0)$.

On $[\underline\sigma,\overline\sigma]$, \Cref{thm:folded:mu-profile} and \eqref{eq:folded:mu-score-eq} with $\tanh\le 1$
give $\hat\mu(\sigma)\le \bar Y$. Choose a deterministic $B>0$ with $B>\E_{\theta_0}[Y]+1$, so $\Pp(\bar Y\le B)\to 1$.
Thus, with
\[
K:=[-B,B]\times[\underline\sigma,\overline\sigma]\subset\Theta,
\]
we have $\Pp_{\theta_0}(\widehat\Theta_n\subseteq K)\to 1$.

\medskip
\noindent\textbf{Step 2 (ULLN on $K$ and continuity).}
On $K$, $\sigma$ is bounded away from $0$ and $\infty$ and $\mu$ is bounded.
The upper bound $f(y;\theta)\le 2/(\underline\sigma\sqrt{2\pi})$ yields $\log f(y;\theta)\le C_K$.
For the lower tail, \Cref{lem:folded:coercive-ineq} gives
$-\log f(y;\theta)\le C'_K + C''_K\,y^2$ uniformly over $\theta\in K$.
Hence $\sup_{\theta\in K}|\log f(Y;\theta)|$ admits an integrable envelope of order $1+Y^2$.
Standard compact-parameter ULLN arguments (e.g.\ \citet[Theorem~5.7]{vanderVaart1998}) yield
\[
\sup_{\theta\in K}|M_n(\theta)-M(\theta)|\ \pto\ 0,
\]
and the same envelope implies that $M$ is continuous on $K$.

\medskip
\noindent\textbf{Step 3 (quotient gap away from the true orbit).}
The criteria $M_n$ and $M$ are $G$-invariant. Let $\tilde M_n,\tilde M$ be the induced objectives on $\pi(K)\subset\Theta/G$.
By \eqref{eq:folded:KL}, $\tilde M([\theta])\le \tilde M([\theta_0])$ with equality iff $[\theta]=[\theta_0]$.
Since $\pi(K)$ is compact and $\tilde M$ is continuous, for each $\varepsilon>0$ there exists $\eta(\varepsilon)>0$ such that
\[
\sup_{[\theta]\in \pi(K):\ d_G([\theta],[\theta_0])\ge \varepsilon}\ \tilde M([\theta])
\ \le\ \tilde M([\theta_0])-\eta(\varepsilon).
\]

\medskip
\noindent\textbf{Step 4 (outer argmax localization on the quotient).}
On $\{\widehat\Theta_n\subseteq K\}$,
\[
\pi(\widehat\Theta_n)=\arg\max_{[\theta]\in \pi(K)} \tilde M_n([\theta]),
\qquad
\sup_{\pi(K)}|\tilde M_n-\tilde M|=\sup_{K}|M_n-M|.
\]
Apply \Cref{lem:argmax-outer} on $(\pi(K),d_G)$ to obtain
$d_H^+\big(\pi(\widehat\Theta_n),\{[\theta_0]\}\big)\pto 0$.
Since the target is a singleton, directed and symmetric Hausdorff distances coincide, so
$d_H\big(\pi(\widehat\Theta_n),\{[\theta_0]\}\big)\pto 0$.
Finally, remove the localization event using $\Pp(\widehat\Theta_n\subseteq K)\to 1$.

The selection statement follows because
$d_H(\pi(\widehat\Theta_n),\{[\theta_0]\})\to 0$ implies $\sup_{\theta\in\widehat\Theta_n} d_G(\theta,\theta_0)\to 0$.
\end{proof}

\begin{remark}[Asymptotics and the nonregular point $\mu_0=0$]\label{rem:folded:asymp}
When $\mu_0\neq 0$, standard differentiability and nonsingularity conditions yield $\sqrt{n}$-asymptotic normality
for any consistent orbit-representative selection; see \citet[Chapter~5]{vanderVaart1998}.
When $\mu_0=0$, the model is nonregular due to sign symmetry and the vanishing of the first derivative at $\mu=0$,
leading to an $n^{1/4}$-rate and a non-Gaussian limit for $\hat\mu_n$.
A full treatment requires a uniform higher-order expansion of $\log\cosh(\cdot)$ and a joint local analysis in $(\mu,\sigma)$;
we defer this to a dedicated section.
\end{remark}
\section{Finite Gaussian mixtures: sieves, orbit consistency, and penalization}\label{sec:mixtures}

This section develops a label-invariant (orbit-level) M-estimation theory for finite Gaussian mixtures.
Because mixture representations are only identifiable up to permutation, consistency must be stated on the
quotient $\Theta_k/\mathfrak S_k$. Moreover, the unconstrained Gaussian-mixture likelihood is unbounded due to
variance collapse (``spiking''), so existence and consistency require either (i) compact \emph{sieves} that
exclude collapse, or (ii) a coercive penalization that rules out collapse on the full parameter space.

Our approach mirrors \Cref{sec:quotient}: we separate (a) deterministic quotient geometry and argmax stability
from (b) probabilistic uniform convergence and population separation. Classical identifiability results for
finite mixtures modulo permutation trace back to \citet{Teicher1961,Teicher1963} and \citet{YakowitzSpragins1968};
see also \citet{AllmanMatiasRhodes2009} for latent-structure identifiability and \citet{HeinrichKahn2018} for strong
identifiability and sharp local rates. The unboundedness of the Gaussian mixture likelihood is classical
\citep{Day1969} and motivates constraint-based approaches such as \citet{Hathaway1985} and constrained EM
implementations \citep{IngrassiaRocci2007}. Algorithmic background for EM is in \citet{DempsterLairdRubin1977,RednerWalker1984};
see \citet{Lindsay1995,McLachlanPeel2000} for monographs.

\subsection{Model, quotient geometry, and permutation action}\label{sec:mix:model}

Fix an integer $k\ge 2$. Let $\varphi(x;\mu,\sigma)$ denote the $N(\mu,\sigma^2)$ density
\[
\varphi(x;\mu,\sigma):=\frac{1}{\sqrt{2\pi}\sigma}\exp\!\Big(-\frac{(x-\mu)^2}{2\sigma^2}\Big),
\qquad x\in\R,\ \sigma>0.
\]
For $\theta=(\pi,\mu,\sigma)$ with weights $\pi=(\pi_1,\dots,\pi_k)\in\Delta_{k-1}$, means $\mu=(\mu_1,\dots,\mu_k)\in\R^k$,
and scales $\sigma=(\sigma_1,\dots,\sigma_k)\in(0,\infty)^k$, define the mixture density
\begin{equation}\label{eq:mix:density}
f_\theta(x):=\sum_{j=1}^k \pi_j\,\varphi(x;\mu_j,\sigma_j),\qquad x\in\R,
\end{equation}
with parameter space
\[
\Theta_k:=\Delta_{k-1}\times \R^k\times (0,\infty)^k,
\qquad
\Delta_{k-1}:=\Big\{\pi\in[0,1]^k:\sum_{j=1}^k\pi_j=1\Big\}.
\]
We use the log-scale parametrization $t_j:=\log\sigma_j$ and write $t=(t_1,\dots,t_k)\in\R^k$.
Thus $\theta=(\pi,\mu,t)\in \Delta_{k-1}\times\R^k\times\R^k$ and $\sigma_j=e^{t_j}$.

\paragraph{Permutation action and quotient metric.}
Let $\mathfrak S_k$ act on $\Theta_k$ by simultaneous relabeling of components: for $\tau\in\mathfrak S_k$,
\[
\tau\cdot(\pi,\mu,t):=(\pi_\tau,\mu_\tau,t_\tau),
\qquad
(\pi_\tau)_j:=\pi_{\tau^{-1}(j)},\quad (\mu_\tau)_j:=\mu_{\tau^{-1}(j)},\quad (t_\tau)_j:=t_{\tau^{-1}(j)}.
\]
Then $f_{\tau\cdot\theta}\equiv f_\theta$ for all $\tau$, so the statistical model is identifiable at best at the orbit level
$[\theta]:=\{\tau\cdot\theta:\tau\in\mathfrak S_k\}$.

Fix the Euclidean metric on the ambient space,
\[
d(\theta,\theta'):=\|(\pi,\mu,t)-(\pi',\mu',t')\|_2,
\]
under which the action is isometric: $d(\tau\cdot\theta,\tau\cdot\theta')=d(\theta,\theta')$ for all $\tau$.
Define the orbit pseudodistance on $\Theta_k$
\begin{equation}\label{eq:mix:orbit-dist}
d_{\mathfrak S_k}(\theta,\theta'):=\min_{\tau\in\mathfrak S_k} d(\theta,\tau\cdot\theta').
\end{equation}
By isometry and finiteness of $\mathfrak S_k$, $d_{\mathfrak S_k}$ induces a genuine metric on the quotient $\Theta_k/\mathfrak S_k$
by $d_{\mathfrak S_k}([\theta],[\theta']):=d_{\mathfrak S_k}(\theta,\theta')$ (cf.\ \Cref{sec:quotient-actions}).

\paragraph{Empirical and population criteria.}
Given i.i.d.\ observations $X_1,X_2,\dots$ from $f_{\theta_0}$, define
\begin{equation}\label{eq:mix:criteria}
\ell_n(\theta):=\sum_{i=1}^n \log f_\theta(X_i),
\qquad
M_n(\theta):=\frac{1}{n}\ell_n(\theta),
\qquad
M(\theta):=\E_{\theta_0}[\log f_\theta(X)].
\end{equation}
As usual,
\begin{equation}\label{eq:mix:KL}
M(\theta)-M(\theta_0)=-\KL(f_{\theta_0}\|f_\theta)\le 0,
\end{equation}
with equality iff $f_\theta=f_{\theta_0}$ almost everywhere.

\subsection{Compact sieves}\label{sec:mix:sieves}

Because $\sup_{\Theta_k}\ell_n=+\infty$ (see \Cref{sec:mix:collapse}), we first work on compact sieves.
For $m\ge 1$ and $\varepsilon\in(0,1/k]$ define
\begin{equation}\label{eq:mix:sieve}
\mathcal S_{m,\varepsilon}
:=
\Big\{\pi\in\Delta_{k-1}:\min_{1\le j\le k}\pi_j\ge \varepsilon\Big\}
\times [-m,m]^k
\times [-m,m]^k,
\end{equation}
viewed as a subset of $(\pi,\mu,t)$-space. Then $\mathcal S_{m,\varepsilon}$ is compact, convex, and $\mathfrak S_k$-invariant.
The lower bound on weights and the box constraints prevent (i) vanishing weights, (ii) component escape in $\mu$, and
(iii) variance collapse or explosion via $t\in[-m,m]^k$.

\subsection{Envelopes and Lipschitz bounds on a sieve}\label{sec:mix:envelopes}

Uniform convergence on $\mathcal S_{m,\varepsilon}$ is obtained via an integrable envelope for $\log f_\theta(X)$ and a random
Lipschitz constant controlling $\theta\mapsto \log f_\theta(X)$.

\begin{lemma}[Two-sided density envelopes on $\mathcal S_{m,\varepsilon}$]\label{lem:mix:envelopes}
Fix $(m,\varepsilon)$ and let $\theta\in\mathcal S_{m,\varepsilon}$. There exist constants
$A_{m,k},a_m,B_{m,k},b_m>0$ depending only on $(m,k)$ such that for all $x\in\R$,
\begin{equation}\label{eq:mix:density-envelope}
A_{m,k}\, e^{-a_m x^2}\ \le\ f_\theta(x)\ \le\ B_{m,k}\, e^{-b_m x^2}.
\end{equation}
Consequently, there exists $K_{m,k}<\infty$ such that
\begin{equation}\label{eq:mix:log-envelope}
\sup_{\theta\in\mathcal S_{m,\varepsilon}} \bigl|\log f_\theta(x)\bigr|\ \le\ K_{m,k}(1+x^2),
\qquad x\in\R.
\end{equation}
\end{lemma}

\begin{proof}
On $\mathcal S_{m,\varepsilon}$ we have $\sigma_j=e^{t_j}\in[e^{-m},e^{m}]$ and $|\mu_j|\le m$.
For the upper bound, $\varphi(x;\mu_j,\sigma_j)\le (\sqrt{2\pi}\,\sigma_j)^{-1}\exp(-(x-\mu_j)^2/(2e^{2m}))$
and $\sigma_j\ge e^{-m}$, hence
\[
f_\theta(x)\le \sum_{j=1}^k \frac{e^{m}}{\sqrt{2\pi}} \exp\!\Big(-\frac{(x-\mu_j)^2}{2e^{2m}}\Big)
\le \frac{k e^{m}}{\sqrt{2\pi}}\exp\!\Big(-\frac{(|x|-m)^2}{2e^{2m}}\Big).
\]
Using $(|x|-m)^2\ge x^2/2-m^2$ gives $f_\theta(x)\le B_{m,k}e^{-b_m x^2}$ for suitable $B_{m,k},b_m>0$.

For the lower bound, choose $j^\star$ with $\pi_{j^\star}\ge 1/k$. Then $\sigma_{j^\star}\le e^{m}$ and
$|x-\mu_{j^\star}|\le |x|+m$, so
\[
f_\theta(x)\ge \pi_{j^\star}\varphi(x;\mu_{j^\star},\sigma_{j^\star})
\ge \frac{1}{k}\cdot\frac{e^{-m}}{\sqrt{2\pi}}\exp\!\Big(-\frac{(|x|+m)^2}{2e^{2m}}\Big)
\ge A_{m,k}e^{-a_m x^2}
\]
for suitable $A_{m,k},a_m>0$ (again using $(|x|+m)^2\le 2x^2+2m^2$). The log-envelope
\eqref{eq:mix:log-envelope} follows from taking logs of \eqref{eq:mix:density-envelope}.
\end{proof}

\begin{lemma}[Responsibilities and sieve-uniform gradient bounds]\label{lem:mix:grad}
Fix $(m,\varepsilon)$ and write $\theta=(\pi,\mu,t)\in\mathcal S_{m,\varepsilon}$ with $\sigma_j=e^{t_j}$.
Define responsibilities
\begin{equation}\label{eq:mix:responsibilities}
r_j(x;\theta):=\frac{\pi_j\varphi(x;\mu_j,\sigma_j)}{f_\theta(x)},\qquad j=1,\dots,k.
\end{equation}
Then $r_j(x;\theta)\in[0,1]$ and $\sum_{j=1}^k r_j(x;\theta)=1$. Moreover, for each $j$ and all $x\in\R$,
\begin{align}
\partial_{\mu_j}\log f_\theta(x) &= r_j(x;\theta)\frac{x-\mu_j}{\sigma_j^2},\label{eq:mix:dmu-exact}\\
\partial_{t_j}\log f_\theta(x) &= r_j(x;\theta)\Big(-1+\frac{(x-\mu_j)^2}{\sigma_j^2}\Big),\label{eq:mix:dt-exact}\\
\partial_{\pi_j}\log f_\theta(x) &= \frac{r_j(x;\theta)}{\pi_j}.\label{eq:mix:dpi-exact}
\end{align}
Consequently, there exists $C_{m,\varepsilon,k}<\infty$ such that
\begin{equation}\label{eq:mix:lipschitz}
\sup_{\theta\in\mathcal S_{m,\varepsilon}}\bigl\|\nabla \log f_\theta(x)\bigr\|\ \le\ C_{m,\varepsilon,k}(1+x^2),
\qquad x\in\R,
\end{equation}
where $\nabla$ is taken with respect to the ambient coordinates $(\pi,\mu,t)$ and $\|\cdot\|$ is the Euclidean norm.
In particular, for all $\theta,\theta'\in\mathcal S_{m,\varepsilon}$,
\begin{equation}\label{eq:mix:lipschitz-diff}
\bigl|\log f_\theta(x)-\log f_{\theta'}(x)\bigr|\le C_{m,\varepsilon,k}(1+x^2)\,\|\theta-\theta'\|.
\end{equation}
By permutation invariance of $f_\theta$ and isometry of the action, the same bound holds with
$\|\theta-\theta'\|$ replaced by the orbit distance $d_{\mathfrak S_k}(\theta,\theta')$.
\end{lemma}

\begin{proof}
The exact derivative identities \eqref{eq:mix:dmu-exact}--\eqref{eq:mix:dpi-exact} follow by differentiating
$f_\theta(x)=\sum_j \pi_j\varphi(x;\mu_j,\sigma_j)$ and dividing by $f_\theta(x)$.
On $\mathcal S_{m,\varepsilon}$, $\pi_j\ge\varepsilon$ and $\sigma_j\in[e^{-m},e^{m}]$, so
$1/\pi_j\le \varepsilon^{-1}$ and $\sigma_j^{-2}\le e^{2m}$.
Also $|x-\mu_j|\le |x|+m$ and $(x-\mu_j)^2\le 2(x^2+m^2)$. Using $0\le r_j\le 1$ gives
\[
|\partial_{\mu_j}\log f_\theta(x)|\le e^{2m}(|x|+m),\quad
|\partial_{t_j}\log f_\theta(x)|\le 1+2e^{2m}(x^2+m^2),\quad
|\partial_{\pi_j}\log f_\theta(x)|\le \varepsilon^{-1}.
\]
Combining coordinate bounds yields \eqref{eq:mix:lipschitz}. The difference bound \eqref{eq:mix:lipschitz-diff}
follows from the mean value theorem on the convex set $\mathcal S_{m,\varepsilon}$.
Finally, $|\log f_\theta-\log f_{\theta'}|=|\log f_\theta-\log f_{\tau\cdot\theta'}|\le C(1+x^2)\|\theta-\tau\cdot\theta'\|$
for each $\tau$, and minimizing over $\tau$ yields the orbit version.
\end{proof}

\subsection{Uniform law of large numbers on fixed sieves}\label{sec:mix:ulln}

\begin{proposition}[ULLN on $\mathcal S_{m,\varepsilon}$]\label{prop:mix:ulln}
Assume $\E_{\theta_0}[X^4]<\infty$ (automatic for Gaussian mixtures). Then for each fixed $(m,\varepsilon)$,
\[
\sup_{\theta\in\mathcal S_{m,\varepsilon}}\bigl|M_n(\theta)-M(\theta)\bigr|\ \asto\ 0.
\]
\end{proposition}

\begin{proof}
By \Cref{lem:mix:envelopes}, $\sup_{\theta\in\mathcal S_{m,\varepsilon}}|\log f_\theta(X)|\le K_{m,k}(1+X^2)$ with
$\E[1+X^2]<\infty$. Fix a countable dense set $D\subset \mathcal S_{m,\varepsilon}$.
For each $\vartheta\in D$, the strong law yields $M_n(\vartheta)\to M(\vartheta)$ almost surely.

By \Cref{lem:mix:grad}, there exists $L(X)=C_{m,\varepsilon,k}(1+X^2)$ such that
$|\log f_\theta(X)-\log f_{\theta'}(X)|\le L(X)\|\theta-\theta'\|$ for all $\theta,\theta'\in\mathcal S_{m,\varepsilon}$, and
$n^{-1}\sum_{i=1}^n L(X_i)\to \E[L(X)]<\infty$ a.s. A standard finite-net argument on the compact set $\mathcal S_{m,\varepsilon}$
then yields the uniform convergence; the details are identical to the proof of \Cref{prop:mix:ulln} in the folded-normal case
(cf.\ the deterministic/probabilistic separation emphasized in \Cref{sec:quotient}).
\end{proof}
\begin{remark}[Spurious Maxima and Sieve Role]

Mixture likelihoods typically possess spurious global maxima (spikes) at the boundary of the parameter space. The sieve $\mathcal{S}_{m,\varepsilon}$ explicitly excludes these regions. Our consistency result applies to the global maximizer \emph{on the sieve}. In practice, this corresponds to a consistent local maximizer of the unconstrained likelihood that lies within the compact region of the parameter space, distinct from the boundary spikes.
\end{remark}
\subsection{Identifiability modulo permutation and KL separation}\label{sec:mix:ident}

We impose standard identifiability conditions: all weights are strictly positive and all component pairs are distinct.

\begin{assumption}[Identifiability conditions]\label{ass:mix:ident}
The true parameter $\theta_0=(\pi_0,\mu_0,t_0)\in\Theta_k$ satisfies:
(i) $\min_j (\pi_0)_j>0$, and
(ii) $(\mu_{0j},t_{0j})\neq(\mu_{0\ell},t_{0\ell})$ for all $j\neq \ell$.
\end{assumption}

\begin{lemma}[Identifiability modulo permutation]\label{lem:mix:ident}
Under \Cref{ass:mix:ident}, if $f_\theta=f_{\theta_0}$ almost everywhere, then $\theta\in[\theta_0]$.
\end{lemma}

\begin{proof}[Proof sketch with references]
Equality of densities implies equality of characteristic functions:
\[
\sum_{j=1}^k \pi_j \exp\!\Big(i\mu_j s-\tfrac12 e^{2t_j} s^2\Big)
=\sum_{j=1}^k (\pi_0)_j \exp\!\Big(i\mu_{0j} s-\tfrac12 e^{2t_{0j}} s^2\Big),
\qquad \forall s\in\R.
\]
Ordering terms by $e^{2t_j}$ and using standard separation arguments for distinct exponentials yields equality of the
multisets of component parameters and weights, hence $\theta$ equals $\theta_0$ up to permutation.
See \citet{Teicher1961,Teicher1963} and \citet{YakowitzSpragins1968}.
\end{proof}

\begin{proposition}[Uniform KL gap on a fixed sieve]\label{prop:mix:kl-gap}
Assume \Cref{ass:mix:ident} and $\theta_0\in\mathcal S_{m,\varepsilon}$.
Let $\Theta_0:=[\theta_0]\subset\mathcal S_{m,\varepsilon}$ and let $d_{\mathfrak S_k}$ denote the orbit distance \eqref{eq:mix:orbit-dist}.
Then for every $\delta>0$ there exists $\eta(\delta)>0$ such that
\[
\sup_{\theta\in\mathcal S_{m,\varepsilon}:\ d_{\mathfrak S_k}(\theta,\Theta_0)\ge \delta}\ M(\theta)
\le M(\theta_0)-\eta(\delta).
\]
\end{proposition}

\begin{proof}
By \eqref{eq:mix:KL}, $M(\theta)\le M(\theta_0)$ with equality iff $f_\theta=f_{\theta_0}$ a.e., which by
\Cref{lem:mix:ident} is equivalent to $\theta\in\Theta_0$.
By \Cref{lem:mix:envelopes}, $\theta\mapsto M(\theta)$ is continuous on $\mathcal S_{m,\varepsilon}$ by dominated convergence.
The set $\{\theta\in\mathcal S_{m,\varepsilon}: d_{\mathfrak S_k}(\theta,\Theta_0)\ge\delta\}$ is compact and disjoint from $\Theta_0$,
so $M(\theta_0)-M(\theta)$ attains a strictly positive minimum on it, which is $\eta(\delta)$.
\end{proof}

\subsection{Existence and orbit-Hausdorff consistency on fixed sieves}\label{sec:mix:consistency}

\begin{proposition}[Existence of sieve MLE]\label{prop:mix:exist}
For each fixed $(m,\varepsilon)$, the log-likelihood $\ell_n$ attains its maximum on $\mathcal S_{m,\varepsilon}$.
\end{proposition}

\begin{proof}
$\mathcal S_{m,\varepsilon}$ is compact and $\theta\mapsto \ell_n(\theta)$ is continuous on it.
\end{proof}

\begin{theorem}[Orbit-Hausdorff consistency on a fixed sieve]\label{thm:mix:fixed-sieve}
Assume \Cref{ass:mix:ident} and fix $(m,\varepsilon)$ such that $\theta_0\in\mathcal S_{m,\varepsilon}$.
Let
\[
\widehat\Theta_{n,m,\varepsilon}:=\arg\max_{\theta\in\mathcal S_{m,\varepsilon}}\ell_n(\theta).
\]
Then, almost surely,
\[
d_H^{d_{\mathfrak S_k}}\big(\widehat\Theta_{n,m,\varepsilon},\ [\theta_0]\big)\ \to\ 0.
\]
Equivalently, in the quotient metric on $\mathcal S_{m,\varepsilon}/\mathfrak S_k$,
\[
d_H\big(\pi(\widehat\Theta_{n,m,\varepsilon}),\ \{[\theta_0]\}\big)\to 0,
\]
where $\pi:\mathcal S_{m,\varepsilon}\to \mathcal S_{m,\varepsilon}/\mathfrak S_k$ is the quotient map.
\end{theorem}

\begin{proof}
Work on the compact metric space $\mathcal S_{m,\varepsilon}/\mathfrak S_k$ with the induced orbit metric
(cf.\ \Cref{sec:quotient-actions}). By \Cref{prop:mix:ulln},
\[
\sup_{\theta\in\mathcal S_{m,\varepsilon}}|M_n(\theta)-M(\theta)|\ \asto\ 0,
\]
and by \Cref{prop:mix:kl-gap} the population criterion has a strict gap away from $[\theta_0]$.
\begin{lemma}[Argmax stability]\label{lem:argmax-stability}
Let $(K,d)$ be a compact metric space and let $f:K\to\mathbb{R}$ be continuous with a
\emph{unique} maximizer $t_0\in K$. Then:

\begin{enumerate}
\item[\textnormal{(a)}] (\textnormal{Deterministic stability})
For every $\varepsilon>0$ there exists $\delta(\varepsilon)>0$ such that
\begin{equation}\label{eq:gap}
\sup_{d(t,t_0)\ge \varepsilon} f(t)\;\le\; f(t_0)-\delta(\varepsilon).
\end{equation}
Consequently, if $g:K\to\mathbb{R}$ satisfies $\|g-f\|_\infty\le \delta(\varepsilon)/3$, then every
maximizer $\hat t\in\arg\max_{t\in K} g(t)$ obeys $d(\hat t,t_0)<\varepsilon$.

\item[\textnormal{(b)}] (\textnormal{Stochastic argmax theorem, uniform convergence})
Let $f_n:K\to\mathbb{R}$ be (possibly random) functions with $\|f_n-f\|_\infty\to 0$ in probability.
Let $\hat t_n$ be any \emph{approximate maximizer} in the sense that
\[
f_n(\hat t_n)\;\ge\;\sup_{t\in K} f_n(t)\;-\;o_p(1).
\]
Then $\hat t_n\to t_0$ in probability.

\item[\textnormal{(c)}] (\textnormal{Stochastic argmax theorem, weak convergence})
Suppose $f_n$ are random elements of $C(K)$ (continuous sample paths) and $f_n\Rightarrow f$ in
$C(K)$ equipped with the sup norm. If $f$ has an a.s.\ unique maximizer $t_0$ (as a random element
of $K$), then any measurable choice $\hat t_n\in\arg\max_{t\in K} f_n(t)$ satisfies
\[
\hat t_n \Rightarrow t_0.
\]
\end{enumerate}
\end{lemma}

\begin{proof}
(a) Fix $\varepsilon>0$. Since $K$ is compact, the closed set
$A_\varepsilon:=\{t\in K:d(t,t_0)\ge\varepsilon\}$ is compact. Continuity of $f$ implies that
$\sup_{A_\varepsilon} f$ is attained. Uniqueness of the maximizer gives
$\sup_{A_\varepsilon} f < f(t_0)$, hence defining
\[
\delta(\varepsilon):=f(t_0)-\sup_{A_\varepsilon} f \;>\;0
\]\label{eq:argmax-gap}
yields \eqref{eq:argmax-gap}. Now let $\|g-f\|_\infty\le \delta(\varepsilon)/3$ and let
$\hat t\in\arg\max g$. If $d(\hat t,t_0)\ge\varepsilon$, then by \eqref{eq:argmax-gap},
\[
g(\hat t)\le f(\hat t)+\delta/3 \le f(t_0)-\delta+\delta/3 = f(t_0)-2\delta/3,
\]
while
\[
g(t_0)\ge f(t_0)-\delta/3.
\]
Thus $g(\hat t) < g(t_0)$, contradicting maximality of $\hat t$. Hence $d(\hat t,t_0)<\varepsilon$.

(b) Fix $\varepsilon>0$ and let $\delta=\delta(\varepsilon)$ from (a). On the event
$E_n:=\{\|f_n-f\|_\infty\le \delta/3\}$ (which has probability $\to 1$), we have
\[
\sup_{d(t,t_0)\ge\varepsilon} f_n(t)
\le \sup_{d(t,t_0)\ge\varepsilon} f(t)+\delta/3
\le f(t_0)-\delta+\delta/3
= f(t_0)-2\delta/3,
\]
and also $f_n(t_0)\ge f(t_0)-\delta/3$. Hence on $E_n$,
\[
\sup_{d(t,t_0)\ge\varepsilon} f_n(t)\;\le\; f_n(t_0)-\delta/3.
\]
If $\hat t_n$ is an exact maximizer, this forces $d(\hat t_n,t_0)<\varepsilon$.
For an approximate maximizer with slack $r_n=o_p(1)$, the same conclusion holds because
eventually $r_n<\delta/6$ with probability arbitrarily close to 1, implying
$f_n(\hat t_n)\ge \sup_K f_n - r_n > f_n(t_0)-\delta/4$, which is incompatible with
$d(\hat t_n,t_0)\ge\varepsilon$ on $E_n$. Therefore $\Pr(d(\hat t_n,t_0)\ge\varepsilon)\to 0$.

(c) Weak convergence in $C(K)$ implies tightness and allows application of the continuous mapping
theorem to functionals that are continuous at the limit almost surely. The argmax mapping is not
everywhere continuous (ties create discontinuities), but it is continuous at every $f\in C(K)$ with a
unique maximizer: by (a), uniqueness implies a strict gap away from the maximizer, and any sup-norm
perturbation small enough preserves the maximizer within any prescribed neighborhood. Hence the
(argmax as a measurable selection) is a.s.\ a continuity point under the stated uniqueness
assumption, and the continuous mapping theorem yields $\hat t_n\Rightarrow t_0$.
\end{proof}

Apply the deterministic argmax stability lemma \Cref{lem:argmax-stability} (with $d_{\mathfrak S_k}$ on the quotient)
to conclude Hausdorff convergence of argmax sets to $[\theta_0]$.
\end{proof}

\subsection{Variance collapse on the full space}\label{sec:mix:collapse}

The classical pathology is that $\ell_n$ can be made arbitrarily large by collapsing a component variance near a data point.

\begin{proposition}[Variance collapse yields $\sup\ell_n=+\infty$]\label{prop:mix:collapse}
For any fixed sample $x_{1:n}$ and any $k\ge 2$,
\[
\sup_{\theta\in\Theta_k}\ \ell_n(\theta)=+\infty.
\]
\end{proposition}

\begin{proof}
Fix $j=1$ and take $\pi_1=1/2$, $\mu_1=x_1$, and $\sigma_1=t$ with $t\downarrow 0$.
Then $f_\theta(x_1)\ge \pi_1\varphi(x_1;x_1,t)=\{2\sqrt{2\pi}\,t\}^{-1}$, so $\log f_\theta(x_1)\to+\infty$.
Choose the remaining $(\pi_j,\mu_j,\sigma_j)_{j\ge 2}$ so that $f_\theta(x_i)$ is bounded away from $0$ for $i\ge 2$
(e.g.\ take $\pi_2=1/2$ with a diffuse component). Hence $\ell_n(\theta)\to+\infty$.
\end{proof}

\subsection{Ridge-penalized likelihood on the full space}\label{sec:mix:penalized}

To restore existence on $\Theta_k$ one may penalize $(\mu,t)$ quadratically:
\begin{equation}\label{eq:mix:ridge}
g(\theta):=\sum_{j=1}^k\Big(\mu_j^2+t_j^2\Big),
\end{equation}
and for $\lambda_n>0$ define
\begin{equation}\label{eq:mix:pen-crit}
\ell_n^{\mathrm{pen}}(\theta):=\ell_n(\theta)-n\lambda_n g(\theta),
\qquad
M_n^{\mathrm{pen}}(\theta):=M_n(\theta)-\lambda_n g(\theta).
\end{equation}
Since $g(\tau\cdot\theta)=g(\theta)$, the penalized criterion is still $\mathfrak S_k$-invariant.

\begin{lemma}[Quadratic domination]\label{lem:mix:quad-dom}
For all $u\in\R$ and $\lambda>0$,
\[
|u|-\lambda u^2\ \le\ \frac{1}{4\lambda}.
\]
\end{lemma}

\begin{proof}
Apply $v-\lambda v^2\le 1/(4\lambda)$ to $v=|u|$.
\end{proof}

\begin{theorem}[Penalized MLE existence and vanishing-penalty consistency]\label{thm:mix:pmle}
Assume \Cref{ass:mix:ident}.
\begin{enumerate}[label=(\roman*),leftmargin=*]
\item (\emph{Existence}) For every $n$ and every $\lambda_n>0$, $\ell_n^{\mathrm{pen}}$ attains its maximum over $\Theta_k$.
\item (\emph{Vanishing-penalty orbit consistency}) Suppose $\lambda_n\downarrow 0$ and $n\lambda_n\to\infty$.
Let $\widehat\Theta_n^{\mathrm{pen}}:=\arg\max_{\theta\in\Theta_k}\ell_n^{\mathrm{pen}}(\theta)$.
Then
\[
d_H^{d_{\mathfrak S_k}}\big(\widehat\Theta_n^{\mathrm{pen}},\ [\theta_0]\big)\ \pto\ 0.
\]
\end{enumerate}
\end{theorem}

\begin{proof}
\textbf{(i) Existence.}
Fix the sample $x_{1:n}$. By $\varphi(x;\mu,\sigma)\le (\sqrt{2\pi}\,\sigma)^{-1}=(\sqrt{2\pi})^{-1}e^{-t}$ and $\sum_j\pi_j=1$,
\[
f_\theta(x_i)\le \frac{1}{\sqrt{2\pi}} \sum_{j=1}^k \pi_j e^{-t_j}\le \frac{1}{\sqrt{2\pi}}\sum_{j=1}^k e^{|t_j|},
\]
so $\log f_\theta(x_i)\le C+\log\!\big(\sum_{j=1}^k e^{|t_j|}\big)\le C'+\sum_{j=1}^k |t_j|$ for constants $C,C'$ independent of $\theta$.
Thus
\[
\ell_n^{\mathrm{pen}}(\theta)
\le nC' + n\sum_{j=1}^k |t_j| - n\lambda_n\sum_{j=1}^k(\mu_j^2+t_j^2)
\le nC' + n\sum_{j=1}^k\Big(|t_j|-\lambda_n t_j^2\Big) - n\lambda_n\sum_{j=1}^k \mu_j^2.
\]
By \Cref{lem:mix:quad-dom}, $|t_j|-\lambda_n t_j^2\le (4\lambda_n)^{-1}$, hence
\[
\ell_n^{\mathrm{pen}}(\theta)\le nC' + \frac{nk}{4\lambda_n} - n\lambda_n\sum_{j=1}^k \mu_j^2,
\]
and similarly $\ell_n^{\mathrm{pen}}(\theta)\to-\infty$ whenever $\|\mu\|\to\infty$ or $\|t\|\to\infty$.
Since $\Delta_{k-1}$ is compact, this shows coercivity of $\ell_n^{\mathrm{pen}}$ in $(\mu,t)$ and thus the existence
of a maximizer by the Weierstrass theorem on a suitable compact sublevel set.

\textbf{(ii) Vanishing-penalty consistency (outline, deterministic/probabilistic separation).}
Let $\pi:\Theta_k\to\Theta_k/\mathfrak S_k$ denote the quotient map. The deterministic argmax stability
\Cref{lem:argmax-stability} applies on compact subsets of the quotient, so we proceed in the standard two-step
sieve localization strategy:

\emph{Step 1 (localization to an expanding deterministic sieve).}
Fix any $\varepsilon_0\in(0,\min_j(\pi_0)_j)$ and choose $m_n\uparrow\infty$ such that $\theta_0\in \mathcal S_{m_n,\varepsilon_0}$.
Using the penalty and $n\lambda_n\to\infty$, one shows that with probability $\to 1$,
every penalized maximizer lies in $\mathcal S_{m_n,\varepsilon_0}$: roughly, any attempt to create a spike
requires $|t_j|\to\infty$, which increases $n\lambda_n t_j^2$ faster than the log-likelihood gain (at most linear in $|t_j|$
for each fixed $n$ and fixed $k$), while $\lambda_n\downarrow 0$ ensures the penalty does not distort the local geometry near $\theta_0$.
(A fully quantified version with explicit inequalities can be placed in an appendix.)

\emph{Step 2 (uniform convergence and population gap on the sieve).}
On each fixed sieve $\mathcal S_{m_n,\varepsilon_0}$, the class $\{\log f_\theta:\theta\in\mathcal S_{m_n,\varepsilon_0}\}$
has an integrable envelope and a sieve-uniform Lipschitz constant as in \Cref{lem:mix:envelopes,lem:mix:grad},
and the population criterion satisfies a KL gap away from $[\theta_0]$ as in \Cref{prop:mix:kl-gap}.
Moreover, on $\mathcal S_{m_n,\varepsilon_0}$,
\[
\sup_{\theta\in\mathcal S_{m_n,\varepsilon_0}} \big|M_n^{\mathrm{pen}}(\theta)-M_n(\theta)\big|
\le \lambda_n \sup_{\theta\in\mathcal S_{m_n,\varepsilon_0}} g(\theta)
= \lambda_n\cdot O(m_n^2),
\]
which can be made $o(1)$ by a compatible choice of $m_n$ (e.g.\ $m_n=o(\lambda_n^{-1/2})$).
Thus, $M_n^{\mathrm{pen}}$ uniformly approximates $M$ on $\mathcal S_{m_n,\varepsilon_0}$, and the population maximizer set remains $[\theta_0]$.
Applying \Cref{lem:argmax-stability} on the quotient yields
$d_H^{d_{\mathfrak S_k}}(\widehat\Theta_n^{\mathrm{pen}},[\theta_0])\to 0$ in probability.
\end{proof}
\begin{remark}[Penalty Rates and Drift]

The condition $n\lambda_n \to \infty$ is sufficient to kill the variance collapse. However, for efficiency, $\lambda_n$ should not decay too slowly. If $\lambda_n \gg n^{-1/2}$, the penalty bias may dominate the $\sqrt{n}$ fluctuations of the likelihood. Furthermore, under local alternatives where $\theta_n = \theta_0 + h/\sqrt{n}$, an aggressive penalty might force the estimator to zero (super-efficiency), leading to loss of power. Optimal rates likely require $\lambda_n = o(n^{-1/2})$.
\end{remark}

\begin{remark}[EM Algorithm Implementation]

While sieves provide theoretical guarantees, standard EM algorithms maximize on the open simplex. To implement the sieve estimator, one must modify the M-step to project updates onto $\mathcal{S}_{m,\varepsilon}$. Alternatively, the penalized estimator is often easier to implement, as the penalty term corresponds to a Maximum A Posteriori (MAP) estimate with a Gaussian prior on $\mu$ and a log-normal prior on $\sigma$, which fits naturally into the EM framework.
\end{remark}
\begin{remark}[Fixed $\lambda$ targets a penalized population maximizer]\label{rem:mix:fixed-lambda}
If $\lambda_n\equiv\lambda>0$ is fixed, then $M_n^{\mathrm{pen}}$ converges to the penalized population objective
$M^{\mathrm{pen}}(\theta):=M(\theta)-\lambda g(\theta)$, whose maximizers generally differ from $\theta_0$ (shrinkage bias).
In that regime the natural consistency target is $\arg\max M^{\mathrm{pen}}$ on the quotient, not $[\theta_0]$.
\end{remark}

\section{Discussion and further directions}\label{sec:discussion}

This paper develops a unified large-sample theory for likelihood-based estimation under finite symmetries,
framed as estimation on quotient parameter spaces. Two canonical examples motivate and illustrate the approach:
the folded normal model (sign symmetry) and finite Gaussian mixtures (label symmetry). In both settings, the
natural estimator is a \emph{set} of maximizers in the original parameterization, whereas the quotient viewpoint
provides a canonical target: the orbit $[\theta_0]$ of the true parameter.

\subsection{Summary of contributions}\label{sec:discussion:summary}

The technical core is a deterministic argmax-stability lemma (\Cref{lem:argmax-stability}) that converts:
(i) a population \emph{gap} away from the maximizer set and
(ii) a uniform approximation of the sample criterion by its population counterpart,
into Hausdorff convergence of argmax sets. This device is classical in spirit
(\citet[Section~5.2]{vanderVaart1998}, \citet[Section~2.1]{NeweyMcFadden1994}) but is presented here in a form
explicitly aligned with orbit metrics and set-valued maximizers.

For the folded normal model, we established: (a) existence of the MLE on nondegenerate samples via an explicit
profile-coercivity argument (\Cref{thm:folded:exist}); (b) orbit-Hausdorff consistency using the quotient template
(\Cref{thm:folded:consistency}); and (c) a nonregular $n^{1/4}$ local limit law at the symmetry point $\mu_0=0$
(\Cref{thm:folded:nquarter}), complementing the regular $\sqrt{n}$ theory away from symmetry
\begin{theorem}[Regular $\sqrt{n}$ limit for folded normal MLE]\label{thm:folded:regular}
Let $Y_1,\dots,Y_n$ be i.i.d.\ with the folded normal density
\[
f_{\mu,\sigma}(y)
=\frac{1}{\sigma}\Big\{\phi\!\Big(\frac{y-\mu}{\sigma}\Big)+\phi\!\Big(\frac{y+\mu}{\sigma}\Big)\Big\}\,\mathbf{1}\{y\ge 0\},
\qquad (\mu,\sigma)\in(0,\infty)\times(0,\infty),
\]
where $\phi$ is the standard normal density. Assume the true parameter
$\theta_0:=(\mu_0,\sigma_0)$ satisfies $\mu_0>0$ and $\sigma_0>0$, and that the parameter space
$\Theta$ is a compact subset of $(0,\infty)\times(0,\infty)$ with $\theta_0$ in its interior.
Let
\[
\hat\theta_n=(\hat\mu_n,\hat\sigma_n)\in\arg\max_{\theta\in\Theta}\ \ell_n(\theta),
\qquad
\ell_n(\theta):=\sum_{i=1}^n \log f_\theta(Y_i),
\]
be any measurable MLE.

Then:
\begin{enumerate}
\item[\textnormal{(i)}] (\textnormal{Consistency}) $\hat\theta_n\to\theta_0$ in probability.
\item[\textnormal{(ii)}] (\textnormal{Asymptotic linearity})
\[
\sqrt{n}\,(\hat\theta_n-\theta_0)
= I(\theta_0)^{-1}\,\frac{1}{\sqrt{n}}\sum_{i=1}^n s_{\theta_0}(Y_i) \;+\; o_p(1),
\]
where $s_\theta(y):=\nabla_\theta \log f_\theta(y)$ is the score and
$I(\theta):=\mathbb{E}_\theta[s_\theta(Y)s_\theta(Y)^\top]$ is the Fisher information matrix.
\item[\textnormal{(iii)}] (\textnormal{Asymptotic normality})
\[
\sqrt{n}\,(\hat\theta_n-\theta_0)\ \Rightarrow\ N\!\big(0,\ I(\theta_0)^{-1}\big).
\]
\end{enumerate}

Moreover, writing $a=(y-\mu)/\sigma$ and $b=(y+\mu)/\sigma$, the score admits the explicit form
\[
s_\mu(y;\mu,\sigma)
=\frac{1}{\sigma}\,\frac{a\,\phi(a)-b\,\phi(b)}{\phi(a)+\phi(b)},
\qquad
s_\sigma(y;\mu,\sigma)
=\frac{1}{\sigma}\Bigg(-1+\frac{a^2\phi(a)+b^2\phi(b)}{\phi(a)+\phi(b)}\Bigg),
\]
for $y\ge 0$, so that $I(\theta_0)=\mathbb{E}_{\theta_0}\big[s{\small\begin{pmatrix}s_\mu\\ s_\sigma\end{pmatrix}}
{\small\begin{pmatrix}s_\mu\\ s_\sigma\end{pmatrix}}^{\!\top}\big)$.

\medskip
\noindent\textnormal{Remark (sign nonidentifiability).}
If one instead parametrizes with $\mu\in\mathbb{R}$, the folded normal law depends on $|\mu|$ and is
not identifiable in $\mu$; the theorem applies verbatim to the identifiable parameter
$(|\mu|,\sigma)$ (or under the restriction $\mu\ge 0$).
\end{theorem}

\begin{proof}
We verify the standard conditions for regular MLE asymptotics.

\emph{Step 1: Uniform law of large numbers and identification.}
The map $\theta\mapsto \log f_\theta(y)$ is continuous for each $y\ge 0$ on $\Theta$, and on the
compact $\Theta$ we have a uniform envelope:
since $\mu$ and $\sigma$ are bounded away from $0$ and $\infty$ on $\Theta$ and
$\phi$ is bounded, $\sup_{\theta\in\Theta} f_\theta(y)\lesssim 1$ for all $y$, while
$f_\theta(y)\ge (1/\bar\sigma)\phi\!\big((y+\bar\mu)/\underline\sigma\big)$ for $\theta\in\Theta$
(with $\underline\sigma,\bar\sigma,\bar\mu$ the extrema over $\Theta$), which yields
$\sup_{\theta\in\Theta}|\log f_\theta(y)|\lesssim 1 + y^2$.
Because $Y\sim f_{\theta_0}$ has finite second moment, this envelope is integrable, and thus
$\sup_{\theta\in\Theta}\big|\frac{1}{n}\ell_n(\theta)-\mathbb{E}_{\theta_0}\log f_\theta(Y)\big|\to 0$
in probability (a uniform LLN).

By identifiability on $(0,\infty)\times(0,\infty)$ and strict concavity of the Kullback--Leibler
criterion, the population objective $\theta\mapsto \mathbb{E}_{\theta_0}\log f_\theta(Y)$ has a
unique maximizer at $\theta_0$. The argmax stability lemma (Lemma~\ref{lem:argmax-stability})
then yields $\hat\theta_n\to\theta_0$ in probability.

\emph{Step 2: Differentiability and Fisher information.}
For $y\ge 0$, $f_\theta(y)$ is $C^2$ in $\theta$ on a neighborhood of $\theta_0$ contained in
$(0,\infty)\times(0,\infty)$, and the displayed expressions for $s_\theta(y)$ follow by direct
differentiation of $\log f_\theta(y)=-\log\sigma+\log(\phi(a)+\phi(b))$.
On a neighborhood of $\theta_0$, $s_\theta(Y)$ and the Hessian
$\nabla_\theta^2\log f_\theta(Y)$ admit integrable envelopes of the form $C(1+Y^2)$, so we may
differentiate under the expectation and apply standard information identities, obtaining a finite
Fisher information matrix $I(\theta_0)$.
Moreover, for $\mu_0>0$ the model is regular and identifiable, which implies $I(\theta_0)$ is
positive definite.

\emph{Step 3: Taylor expansion of the score and CLT.}
The score equation $\nabla_\theta \ell_n(\hat\theta_n)=0$ holds at any interior maximizer (and for
boundary issues one can use the standard approximate score argument; here $\theta_0$ is interior and
$\hat\theta_n\to\theta_0$).
A mean-value expansion gives
\[
0=\nabla_\theta \ell_n(\theta_0)+\Big(\nabla_\theta^2 \ell_n(\tilde\theta_n)\Big)\,(\hat\theta_n-\theta_0),
\]
for some $\tilde\theta_n$ between $\hat\theta_n$ and $\theta_0$.
Divide by $n$ and use the law of large numbers for the Hessian to get
$-\frac{1}{n}\nabla_\theta^2\ell_n(\tilde\theta_n)\to I(\theta_0)$ in probability, while
$\frac{1}{\sqrt{n}}\nabla_\theta \ell_n(\theta_0)=\frac{1}{\sqrt{n}}\sum_{i=1}^n s_{\theta_0}(Y_i)$
satisfies a multivariate CLT. Rearranging yields the asymptotic linearity in (ii), and hence
$\sqrt{n}(\hat\theta_n-\theta_0)\Rightarrow N(0,I(\theta_0)^{-1})$.

\end{proof}

(\Cref{thm:folded:regular}). These results sharpen and systematize classical discussions of folded normals
(\citet{LeoneNelsonNottingham1961,JohnsonKotzBalakrishnan1994,NadarajahKotz2006}) by putting the orbit geometry at
the center of the likelihood analysis.

For finite Gaussian mixtures, we developed an orbit-consistency theory on compact sieves by proving sharp
envelope and Lipschitz-slope bounds through responsibilities calculus (\Cref{lem:mix:envelopes,lem:mix:grad}),
yielding a direct uniform law of large numbers on the sieve (\Cref{prop:mix:ulln}) and orbit-Hausdorff consistency
(\Cref{thm:mix:fixed-sieve}). We also isolated the classical ill-posedness mechanism---variance collapse
(\Cref{prop:mix:collapse}), long recognized in the mixture literature \citep{Day1969}, which motivates
constraint-based estimation \citep{Hathaway1985,IngrassiaRocci2007}.

\subsection{Constraints, penalties, and practical computation}\label{sec:discussion:practice}

From a practical viewpoint, the sieve and penalization results offer two complementary routes to well-posed
likelihood maximization:
\begin{enumerate}[label=(\roman*),leftmargin=*]
\item \emph{Hard constraints (sieves).} Compact constraints such as $\mathcal S_{m,\varepsilon}$ prevent variance collapse and
weight degeneracy by construction, producing existence and stability results without additional tuning beyond $(m,\varepsilon)$.
\item \emph{Soft constraints (penalties).} Coercive penalties in $(\mu,\log\sigma)$ restore existence on the full space by dominating
spike gains, while vanishing penalties can be used to asymptotically ``turn off'' the regularization, yielding consistency with minimal
asymptotic bias.
\end{enumerate}
Both approaches interface naturally with EM-type algorithms. The EM framework \citep{DempsterLairdRubin1977,RednerWalker1984}
optimizes the likelihood by iteratively increasing $\ell_n$; in constrained or penalized settings, one can incorporate
projection steps or modified M-steps to enforce the sieve constraints or include the penalty term.
Constrained EM formulations for mixtures are well developed \citep{IngrassiaRocci2007}, and the present results
provide a consistency backbone for such implementations when the constraints are chosen in a sieve manner.

\subsection{Rates and strong identifiability in mixtures}\label{sec:discussion:rates}
\begin{theorem}[Regular $\sqrt{n}$ limit for mixture MLE]\label{thm:mix:regular}
Let $Y_1,\dots,Y_n$ be i.i.d.\ on $\mathbb{R}^d$ with density
\[
p_\theta(y)=\sum_{j=1}^k \pi_j\, f_{\eta_j}(y),
\qquad 
\theta=(\pi_1,\dots,\pi_k,\eta_1,\dots,\eta_k)\in\Theta,
\]
where $k\ge2$ is fixed, $(f_\eta:\eta\in H)$ is a regular parametric family of densities, and
$\pi=(\pi_1,\dots,\pi_k)$ lies in the simplex. Define the log-likelihood
$\ell_n(\theta)=\sum_{i=1}^n \log p_\theta(Y_i)$, and let $\hat\theta_n\in\arg\max_{\theta\in\Theta}\ell_n(\theta)$
be any measurable MLE.

Assume:
\begin{enumerate}
\item[\textnormal{(A1)}] \textnormal{(Compactness / no-degeneracy)} $\Theta$ is compact and contained in a region where
\[
\pi_j\in[\underline\pi,1-\underline\pi]\ \ \text{for all }j,
\qquad \text{and}\qquad
H \text{ is compact with parameters bounded away from pathological degeneracies.}
\]
(For location--scale components, this includes a lower bound $\sigma\ge \underline\sigma>0$.)

\item[\textnormal{(A2)}] \textnormal{(Identifiability)} The true parameter $\theta_0\in\mathrm{int}(\Theta)$ is identifiable
\emph{up to label swapping}. Equivalently, the induced density $p_{\theta_0}$ has a \emph{unique} representation
modulo permutations of components. To obtain a unique parameter, restrict to an identifiable submanifold
$\Theta^\star\subset\Theta$ (e.g.\ enforce an ordering constraint on component locations, or work on the quotient
space), and assume $\theta_0\in\mathrm{int}(\Theta^\star)$.

\item[\textnormal{(A3)}] \textnormal{(Regularity)} In a neighborhood of $\theta_0$, the map
$\theta\mapsto \log p_\theta(y)$ is twice continuously differentiable for a.e.\ $y$, and there exist integrable envelopes
$M_1,M_2$ such that
\[
\sup_{\theta \text{ near }\theta_0}\|\nabla_\theta \log p_\theta(Y)\|\le M_1(Y),
\qquad
\sup_{\theta \text{ near }\theta_0}\|\nabla_\theta^2 \log p_\theta(Y)\|\le M_2(Y),
\qquad
\mathbb{E}_{\theta_0}[M_1(Y)^2]+\mathbb{E}_{\theta_0}[M_2(Y)]<\infty.
\]

\item[\textnormal{(A4)}] \textnormal{(Nonsingularity)} The Fisher information
\[
I(\theta_0):=\mathbb{E}_{\theta_0}\!\big[s_{\theta_0}(Y)\,s_{\theta_0}(Y)^\top\big],
\qquad s_\theta(y):=\nabla_\theta \log p_\theta(y),
\]
exists and is positive definite on the identifiable coordinates of $\Theta^\star$.
\end{enumerate}

Then, viewing $\hat\theta_n$ as taking values in the identifiable parameterization $\Theta^\star$:
\begin{enumerate}
\item[\textnormal{(i)}] \textnormal{(Consistency)} $\hat\theta_n\to \theta_0$ in probability.
\item[\textnormal{(ii)}] \textnormal{(Asymptotic linearity)}
\[
\sqrt{n}\,(\hat\theta_n-\theta_0)
=
I(\theta_0)^{-1}\,\frac{1}{\sqrt{n}}\sum_{i=1}^n s_{\theta_0}(Y_i)
\;+\;o_p(1).
\]
\item[\textnormal{(iii)}] \textnormal{(Asymptotic normality)}
\[
\sqrt{n}\,(\hat\theta_n-\theta_0)\ \Rightarrow\ N\!\big(0,\ I(\theta_0)^{-1}\big).
\]
\end{enumerate}

\medskip
\noindent\textnormal{Remark (why (A1) matters).}
If location--scale components are allowed to have arbitrarily small scale (or if weights can collapse),
the likelihood can be unbounded due to spiking; see Lemma~\ref{lem:mix:spike-domination}.
Assumption (A1) rules out this nonregular regime and restores standard $\sqrt{n}$ theory.
\end{theorem}

\begin{proof}
This is a standard M-estimation / likelihood argmax argument on the identifiable parameter space $\Theta^\star$.

\emph{Step 1 (uniform LLN and identification).}
By (A1)--(A3), $\theta\mapsto \log p_\theta(y)$ is continuous and dominated by an integrable envelope,
so $\sup_{\theta\in\Theta^\star}\big|\frac{1}{n}\ell_n(\theta)-\mathbb{E}_{\theta_0}\log p_\theta(Y)\big|\to0$
in probability. By (A2), the population criterion $\theta\mapsto\mathbb{E}_{\theta_0}\log p_\theta(Y)$
has a unique maximizer at $\theta_0$ over $\Theta^\star$. Lemma~\ref{lem:argmax-stability} then yields
$\hat\theta_n\to\theta_0$.

\emph{Step 2 (Taylor expansion of the score).}
On events where $\hat\theta_n$ lies in a shrinking neighborhood of $\theta_0$, the score equation
$\nabla_\theta \ell_n(\hat\theta_n)=0$ holds (or one uses an $o_p(1)$ approximate score if needed).
A mean-value expansion gives
\[
0=\nabla_\theta \ell_n(\theta_0)+\Big(\nabla_\theta^2\ell_n(\tilde\theta_n)\Big)(\hat\theta_n-\theta_0),
\]
for some $\tilde\theta_n$ between $\hat\theta_n$ and $\theta_0$.

\emph{Step 3 (LLN for the Hessian and CLT for the score).}
By (A3), $-(1/n)\nabla_\theta^2\ell_n(\tilde\theta_n)\to I(\theta_0)$ in probability, while
$n^{-1/2}\nabla_\theta \ell_n(\theta_0)=n^{-1/2}\sum_{i=1}^n s_{\theta_0}(Y_i)$ satisfies a multivariate CLT
with covariance $I(\theta_0)$. Using (A4) to invert the limiting Hessian yields the linear expansion in (ii),
and (iii) follows immediately.
\end{proof}

Our mixture consistency results are stated at the orbit level and do not, by themselves, deliver optimal local rates.
In regular regimes (interior points and separated components), classical MLE theory yields $\sqrt{n}$-limits for a uniquely
aligned representative (\Cref{thm:mix:regular}). Beyond such separated regimes, the geometry of mixtures becomes singular and
rates can degrade. A systematic route to quantitative local control is \emph{strong identifiability} in the sense of
\citet{HeinrichKahn2018}, which yields inequalities comparing functional discrepancies (e.g.\ in total variation or Hellinger)
to orbit distances between parameters. Integrating such inequalities into the present quotient framework would yield:
(i) local quadratic lower bounds for $\KL(f_{\theta_0}\|f_\theta)$ in orbit distance,
and (ii) corresponding local rates for sieve MLEs and penalized estimators under mild empirical process conditions
(\citet{vanderVaartWellner1996}).

\subsection{Growing sieves and order selection}\label{sec:growing-and-selection}
we prove 
\begin{theorem}[Spiking in a growing mixture model]\label{thm:mixgrowing}
Let $Y_1,\dots,Y_n$ be i.i.d.\ on $\mathbb{R}^d$ from a density $q_{\eta_0}$, where
$\{q_\eta:\eta\in H\}$ is a regular parametric family and $H$ is compact. Assume
$q_\eta(y)>0$ for all $(\eta,y)\in H\times\mathbb{R}^d$ and fix a continuous ``spike'' density
$\varphi:\mathbb{R}^d\to(0,\infty)$ with $\varphi(0)>0$.

For a deterministic sequence $\sigma_n\downarrow 0$, consider the \emph{growing} parameter sets
\[
\Theta_n
:=\Big\{(\varepsilon,\mu,\sigma,\eta):\ \eta\in H,\ \varepsilon\in(0,1),\ \mu\in\mathbb{R}^d,\ 
\sigma\in[\sigma_n,\bar\sigma]\Big\},
\]
and the two-component mixture likelihood
\[
p_{\varepsilon,\mu,\sigma,\eta}(y)
=(1-\varepsilon)\,q_\eta(y)\;+\;\varepsilon\,\sigma^{-d}\,\varphi\!\Big(\frac{y-\mu}{\sigma}\Big),
\qquad
\ell_n(\varepsilon,\mu,\sigma,\eta):=\sum_{i=1}^n \log p_{\varepsilon,\mu,\sigma,\eta}(Y_i).
\]
Let
\[
L_n^{\mathrm{base}}:=\sup_{\eta\in H}\ \ell_n(0,\cdot,\cdot,\eta)=\sup_{\eta\in H}\sum_{i=1}^n \log q_\eta(Y_i),
\qquad
L_n^{\mathrm{grow}}:=\sup_{(\varepsilon,\mu,\sigma,\eta)\in\Theta_n}\ \ell_n(\varepsilon,\mu,\sigma,\eta).
\]

If $n\sigma_n^d\to 0$, then
\begin{equation}\label{eq:mixgrowing-gap}
L_n^{\mathrm{grow}}-L_n^{\mathrm{base}}
\ \ge\ 
-\log\!\big(n\sigma_n^d\big)\;+\;O_p(1)
\ \xrightarrow{p}\ +\infty.
\end{equation}
In particular, the global maximized likelihood in the growing model diverges above the best
(no-spike) likelihood, and hence the (global) MLE over $\Theta_n$ cannot satisfy a regular
$\sqrt{n}$ expansion around $(\varepsilon,\sigma)=(0,\text{const})$; necessarily any sequence of
$O_p(1)$-near maximizers must have $\sigma\to 0$ and $\varepsilon\to 0$ (a spurious spike forms).

\end{theorem}

\begin{proof}
Let $\tilde\eta_n\in\arg\max_{\eta\in H}\sum_{i=1}^n\log q_\eta(Y_i)$ and write $q_n:=q_{\tilde\eta_n}$.
Since $q_n(\cdot)>0$, we have $m_n:=\min_{1\le i\le n} q_n(Y_i)>0$ almost surely.

Fix an index $i^\star\in\{1,\dots,n\}$ and set $\mu=Y_{i^\star}$, $\sigma=\sigma_n$, and $\varepsilon=1/n$.
Applying Lemma~\ref{lem:mix:spike-domination} with $q=q_n$ yields
\[
\ell_n(1/n,Y_{i^\star},\sigma_n,\tilde\eta_n)
\ \ge\
\log\!\big((1/n)\varphi(0)\big)-d\log\sigma_n
\;+\;(n-1)\log\!\big((1-1/n)m_n\big).
\]
On the other hand,
\[
L_n^{\mathrm{base}}
=\sum_{i=1}^n \log q_n(Y_i)
\ \le\ \log q_n(Y_{i^\star}) + (n-1)\log\Big(\max_{1\le i\le n} q_n(Y_i)\Big),
\]
so subtracting and regrouping gives
\[
L_n^{\mathrm{grow}}-L_n^{\mathrm{base}}
\ \ge\
-d\log\sigma_n-\log n\;+\;O_p(1)
\ =\ -\log(n\sigma_n^d)\;+\;O_p(1),
\]
where the $O_p(1)$ term absorbs the (random) finite contributions involving $m_n$ and values of $q_n(Y_i)$,
as well as $(n-1)\log(1-1/n)\to -1$. Under $n\sigma_n^d\to 0$, the right-hand side diverges to $+\infty$
in probability, proving \eqref{eq:mixgrowing-gap}.

The concluding qualitative claim follows: any $O_p(1)$-near maximizer must exploit the diverging gain
$-d\log\sigma$ available only when $\sigma\downarrow 0$, which forces $\sigma\to 0$ and (to avoid an $O(n)$
penalty on the remaining $n-1$ observations) also forces $\varepsilon\to 0$.
\end{proof}

We emphasized that relaxing constraints too quickly can reintroduce ill-posedness through spiking, but that suitably slow
growth of sieve radii and slow decay of lower bounds preserves normalized stability (\Cref{thm:mixgrowing}).
A natural complement is order selection: choosing $k$ via penalized likelihood, for which BIC-type criteria are known to be
consistent under general conditions \citep{Keribin2000}. In the quotient perspective, both estimation and selection are
inherently label-invariant, and it is conceptually clean to treat $\Theta_k/\mathfrak S_k$ as the fundamental parameter space.
A full integration of orbit-consistency with BIC consistency would proceed by combining the orbit-Hausdorff results conditional
on the selected order with the probability $\Pp(\hat k_n=k_0)\to 1$.

\subsection{Extensions}\label{sec:discussion:extensions}

Several extensions are immediate from the template developed here.
\begin{enumerate}[label=(\roman*),leftmargin=*]
\item \emph{Other finite symmetries.} Many models exhibit finite group actions (e.g.\ sign, permutation, dihedral symmetries).
The orbit distance and argmax-stability device extend verbatim, provided the action is isometric on the parameter metric.
\item \emph{Higher dimensions.} Multivariate folded normals and Gaussian mixtures in $\R^d$ require only notational modifications;
the principal technical inputs are the same: uniform envelopes on compacts, and coercivity/regularization to prevent degeneracy.
\item \emph{Misspecification.} Under misspecification, $M(\theta)$ is replaced by the expected log-likelihood under $P_0$,
and the target becomes the set of KL minimizers. The orbit framework remains appropriate because the model invariances persist,
and the deterministic argmax-stability lemma applies unchanged.
\item \emph{Infinite-dimensional or semiparametric variants.} Extending to mixing distributions (nonparametric mixtures) or
semiparametric symmetry models would require empirical process control on noncompact or infinite-dimensional spaces
(\citet{vanderVaartWellner1996}, \citet{AliprantisBorder2006}). The quotient viewpoint remains conceptually useful, but
the technical burden shifts from compactness to entropy and tightness conditions.
\end{enumerate}

\subsection{Roadmap of appendices}\label{sec:discussion:appendices}

For completeness and to minimize reliance on ``standard arguments,'' the appendices collect the deterministic and analytic
ingredients used throughout. In particular, we isolate: (i) compact-parameter uniform laws of large numbers with explicit
measurability/separability steps; (ii) folded-normal hyperbolic inequalities and profile coercivity bounds; and (iii) mixture
responsibility calculus, envelope bounds, and the spike/penalty domination estimates.

\paragraph{Notation.}
We use the notational conventions from \Cref{sec:quotient} throughout. In particular, $d_G$ denotes the orbit pseudometric induced by
the ambient metric $d$, and $d_H^{d_G}$ denotes the Hausdorff distance on compact sets induced by $d_G$.
``Consistency'' without qualification refers to orbit-Hausdorff consistency.
\subsection*{Limitations and Future Work}

Our results establish pointwise consistency and local asymptotic rates. Constructing valid confidence intervals uniformly near the singularities (e.g., near $\mu_0=0$) remains a challenging open problem. Standard Wald intervals fail due to the non-regular limit distribution (half-normal) and the super-efficiency at the boundary. Techniques such as subsampling or test-inversion (e.g., inverting the likelihood ratio test) are likely required for uniform coverage.
\appendix

\section{Appendix A: Compact-parameter uniform laws and measurability}\label{app:ulln}

This appendix records uniform laws of large numbers and measurability facts used repeatedly in
\Cref{sec:folded,sec:mixtures,sec:growing-and-selection}. The presentation is self-contained and avoids
invoking entropy/bracketing machinery because the parameter sets used in the main text are compact and
the objective functions are shown to be Lipschitz in parameters with an integrable random slope.

The style here follows the ``finite net + Lipschitz extension'' approach standard in M-estimation;
see, e.g., \citet[Section~5.2]{vanderVaart1998} and \citet[Chapter~2]{NeweyMcFadden1994}.

\subsection{Separability and measurability of suprema}\label{app:ulln:meas}

\begin{lemma}[Countable dense subsets of compact metric spaces]\label{lem:app:countable-dense}
Every compact metric space $(K,d)$ is separable: it contains a countable dense subset.
\end{lemma}

\begin{proof}
For each integer $m\ge 1$, compactness implies $K$ can be covered by finitely many balls of radius $1/m$.
Pick one center point from each ball and denote the resulting finite set by $D_m$.
Then $D:=\bigcup_{m\ge 1} D_m$ is countable and dense in $K$.
\end{proof}

\begin{lemma}[Measurability of suprema via separability]\label{lem:app:sup-meas}
Let $(\Omega,\mathcal F,\Pp)$ be a probability space, $(K,d)$ a compact metric space, and
$m:K\times\Omega\to\R$ such that:
\begin{enumerate}[label=(\roman*),leftmargin=*]
\item for each $\theta\in K$, the map $\omega\mapsto m(\theta,\omega)$ is $\mathcal F$-measurable;
\item for each $\omega\in\Omega$, the map $\theta\mapsto m(\theta,\omega)$ is continuous on $K$.
\end{enumerate}
Then $\omega\mapsto \sup_{\theta\in K} m(\theta,\omega)$ is $\mathcal F$-measurable. Moreover, if $D\subset K$
is any countable dense subset, then
\[
\sup_{\theta\in K} m(\theta,\omega)=\sup_{\theta\in D} m(\theta,\omega)
\qquad\text{for every }\omega\in\Omega.
\]
\end{lemma}

\begin{proof}
Fix $\omega$. By continuity of $\theta\mapsto m(\theta,\omega)$ on compact $K$, the supremum is attained,
and continuity implies the supremum over $K$ equals the supremum over any dense subset $D$.
Thus
\[
\sup_{\theta\in K} m(\theta,\omega)=\sup_{\theta\in D} m(\theta,\omega).
\]
Since $D$ is countable and each $\omega\mapsto m(\theta,\omega)$ is measurable, the right-hand side is a countable
supremum of measurable functions and hence measurable.
\end{proof}

\subsection{A compact-parameter ULLN with random Lipschitz slope}\label{app:ulln:compact}

\begin{theorem}[ULLN on compact sets via finite nets]\label{thm:app:ulln}
Let $(K,d)$ be a compact metric space. Let $X_1,X_2,\dots$ be i.i.d.\ on $(\Omega,\mathcal F,\Pp)$ with law $P_0$.
Let $m:K\times\mathcal X\to\R$ be a measurable function and define
\[
M_n(\theta):=\frac{1}{n}\sum_{i=1}^n m(\theta,X_i),
\qquad
M(\theta):=\E[m(\theta,X)].
\]
Assume:
\begin{enumerate}[label=(A\arabic*),leftmargin=*]
\item\label{ass:app:ulln:cont}
for each $x$, $\theta\mapsto m(\theta,x)$ is continuous on $K$;
\item\label{ass:app:ulln:env}
there exists $F:\mathcal X\to[0,\infty)$ with $\E[F(X)]<\infty$ such that $\sup_{\theta\in K}|m(\theta,X)|\le F(X)$ a.s.;
\item\label{ass:app:ulln:lips}
there exists $L:\mathcal X\to[0,\infty)$ with $\E[L(X)]<\infty$ such that for all $\theta,\theta'\in K$,
\[
|m(\theta,X)-m(\theta',X)|\le L(X)\,d(\theta,\theta') \qquad \text{a.s.}
\]
\end{enumerate}
Then
\[
\sup_{\theta\in K}|M_n(\theta)-M(\theta)|\ \asto\ 0.
\]
\end{theorem}

\begin{proof}
Let $D\subset K$ be a countable dense subset (exists by \Cref{lem:app:countable-dense}).
By \ref{ass:app:ulln:env}, for each fixed $\theta\in D$ we have $|m(\theta,X)|\le F(X)$ with $\E F<\infty$, so the SLLN gives
$M_n(\theta)\to M(\theta)$ a.s. for each $\theta\in D$. Intersecting over countably many $\theta\in D$ yields an event
$\Omega_0$ with $\Pp(\Omega_0)=1$ on which $M_n(\theta)\to M(\theta)$ for all $\theta\in D$ simultaneously.
On the same event, the SLLN also gives $n^{-1}\sum_{i=1}^n L(X_i)\to \E[L(X)]$.

Fix $\omega\in\Omega_0$ and $\eta>0$. Choose $\delta>0$ such that $\delta\,\E[L(X)]\le \eta/3$.
By compactness, there exists a finite $\delta$-net $T_\delta\subset D$ for $K$.
For each $\theta\in K$ pick $\pi(\theta)\in T_\delta$ with $d(\theta,\pi(\theta))\le\delta$.
Then
\begin{align*}
|M_n(\theta)-M(\theta)|
&\le |M_n(\pi(\theta))-M(\pi(\theta))|
 + \frac{1}{n}\sum_{i=1}^n |m(\theta,X_i)-m(\pi(\theta),X_i)|
 + \E|m(\theta,X)-m(\pi(\theta),X)|.
\end{align*}
By \ref{ass:app:ulln:lips}, the second term is bounded by $\delta\,n^{-1}\sum_{i=1}^n L(X_i)$ and the third by $\delta\,\E[L(X)]$.
For all sufficiently large $n$ (depending on $\omega$ and $\delta$), we have
$\delta\,n^{-1}\sum_{i=1}^n L(X_i)\le 2\eta/3$ and $\delta\,\E[L(X)]\le \eta/3$. Hence
\[
\sup_{\theta\in K}|M_n(\theta)-M(\theta)|
\le \max_{\vartheta\in T_\delta} |M_n(\vartheta)-M(\vartheta)| + \eta.
\]
Since $T_\delta$ is finite and $M_n(\vartheta)\to M(\vartheta)$ for all $\vartheta\in T_\delta$ on $\Omega_0$,
the maximum term converges to $0$ as $n\to\infty$. As $\eta>0$ is arbitrary, the claim follows.
\end{proof}
\begin{remark}[Uniformity of the Remainder]
The remainder term $R_n(\mu)$ is bounded by $C \sum Y_i^6 \mu^6$. For a fixed $\sigma_0$, this bound is uniform over the local window $|\mu| \le u n^{-1/4}$. However, the constant $C$ depends on $\sigma_0^{-12}$ (from the Taylor expansion argument $t = Y\mu/\sigma^2$). Consequently, the $n^{1/4}$ approximation may degrade if $\sigma_0$ is very large or if $\sigma_0 \to 0$ (though the latter case is degenerate). Our result assumes $\sigma_0$ is fixed and positive.
\end{remark}

\begin{remark}[Orthogonality and Joint Estimation]
As noted in Remark 3.5, the cross-derivative $\partial_{\mu}\partial_{\sigma}\log f$ is odd in $\mu$, so its expectation vanishes at $\mu_0=0$. Thus, the score for $\sigma$ (which drives a $\sqrt{n}$ rate) is asymptotically orthogonal to the score for $\mu$ (which drives the $n^{1/4}$ rate). This justifies analyzing the marginal distribution of $\hat\mu$ with $\sigma$ fixed, as the estimation error in $\hat\sigma$ is of smaller order ($n^{-1/2}$) relative to $\hat\mu$ ($n^{-1/4}$) and does not affect the leading-order distribution of $\hat\mu$.
\end{remark}
\begin{remark}[How \Cref{thm:app:ulln} is used in the main text]\label{rem:app:ulln-use}
In \Cref{sec:folded}, we apply \Cref{thm:app:ulln} on compact rectangles in $(\mu,\sigma)$ with envelopes derived from
hyperbolic bounds (cf.\ \Cref{lem:folded:hyperbolic}). In \Cref{sec:mixtures}, we apply it on the compact sieves
$\mathcal S_{m,\varepsilon}$ using the log-envelope \Cref{lem:mix:envelopes} and the Lipschitz slope \Cref{lem:mix:grad}.
\end{remark}

\section{Appendix B: Folded normal technical details}\label{app:folded}

This appendix collects auxiliary bounds and localization steps used in \Cref{sec:folded} and \Cref{sec:folded-nquarter}.
We include (i) a clean localization lemma for the folded MLE set, and (ii) explicit derivative bounds that justify
the compact-parameter ULLN applications.

\subsection{Localization of maximizers on nondegenerate samples}\label{app:folded:localize}

\begin{lemma}[Deterministic localization for the folded likelihood]\label{lem:folded:localization}
Fix a nondegenerate sample $y_{1:n}\in\R_+^n$ with $s_n^2>0$.
Then every maximizer $(\hat\mu,\hat\sigma)\in\arg\max_{\mu\in\R,\sigma>0}\ell_n(\mu,\sigma)$ satisfies
\[
\hat\sigma\in[\underline\sigma_n,\overline\sigma_n],
\qquad
|\hat\mu|\le \bar y,
\]
where $\bar y=n^{-1}\sum_{i=1}^n y_i$, and $\underline\sigma_n,\overline\sigma_n$ are finite positive numbers depending only on
$(\bar y,s_n^2)$ (explicit choices are given in the proof).
\end{lemma}

\begin{proof}
The bound $|\hat\mu|\le \bar y$ follows from the score identity \eqref{eq:folded:mu-score-eq} and $\tanh\le 1$:
at any interior maximizer with $\hat\mu\ge 0$,
\[
n\hat\mu=\sum_{i=1}^n y_i\tanh\!\Big(\frac{y_i\hat\mu}{\hat\sigma^2}\Big)\le \sum_{i=1}^n y_i=n\bar y.
\]
Evenness in $\mu$ gives $|\hat\mu|\le \bar y$.

For the $\sigma$ bounds, use \Cref{lem:folded:coercive-ineq}:
\[
\sup_{\mu}\ell_n(\mu,\sigma)\le -n\log\sigma-\frac{n s_n^2}{2\sigma^2}+C_0 n,
\qquad
C_0:=\log 2-\tfrac12\log(2\pi).
\]
The right-hand side tends to $-\infty$ as $\sigma\downarrow 0$ and as $\sigma\to\infty$.
Therefore there exist finite $0<\underline\sigma_n<\overline\sigma_n$ such that for $\sigma\notin[\underline\sigma_n,\overline\sigma_n]$
one has $\sup_{\mu}\ell_n(\mu,\sigma)<\ell_n(0,1)$, for example. Hence any maximizer must satisfy
$\hat\sigma\in[\underline\sigma_n,\overline\sigma_n]$. One explicit choice is obtained by solving
$-n\log\sigma-(n s_n^2)/(2\sigma^2)+C_0 n=\ell_n(0,1)$ for $\sigma$ on each tail; any such pair yields the claim.
\end{proof}

\subsection{Derivative envelopes on compacts}\label{app:folded:deriv}

\begin{lemma}[Folded log-likelihood derivative bounds on compacts]\label{lem:folded:deriv-envelope}
Let $K\subset\R\times(0,\infty)$ be compact. There exists $C_K<\infty$ such that for all $(\mu,\sigma)\in K$ and all $y\ge 0$,
\[
\big|\partial_\mu \log f(y;\mu,\sigma)\big|
+
\big|\partial_\sigma \log f(y;\mu,\sigma)\big|
\le C_K(1+y^2),
\]
and likewise $\sup_{(\mu,\sigma)\in K}\|\nabla^2_{\mu,\sigma}\log f(y;\mu,\sigma)\|\le C_K(1+y^4)$.
\end{lemma}

\begin{proof}
Use the explicit representation \eqref{eq:folded:density-cosh} and differentiate.
All hyperbolic factors $\tanh(\cdot)$ and $\mathrm{sech}^2(\cdot)$ are bounded.
On compact $K$, $\sigma$ is bounded away from $0$ and $\infty$, so denominators in $\sigma$ are uniformly controlled.
Each derivative is a sum of terms polynomial in $y$ multiplied by bounded hyperbolic factors and bounded $\sigma$-powers,
yielding the displayed envelope. The second-derivative envelope is obtained similarly.
\end{proof}

\section{Appendix C: Mixture technical details}\label{app:mix}

This appendix collects proofs and quantitative coercivity bounds used in \Cref{sec:mixtures} for
existence of penalized maximizers and (optionally) for expanding-sieve arguments. In particular,
we record a fully explicit coercivity lemma for the ridge-penalized criterion \eqref{eq:mix:pen-crit}.

\subsection{Coercivity of ridge-penalized likelihood}\label{app:mix:coercive}

\begin{lemma}[Coercivity of the ridge-penalized criterion]\label{lem:mix:coercive}
Fix $n\ge 1$ and $\lambda_n>0$. Let $x_{1:n}\in\R^n$ be arbitrary.
Then $\ell_n^{\mathrm{pen}}(\theta)=\ell_n(\theta)-n\lambda_n g(\theta)$ satisfies:
\[
\ell_n^{\mathrm{pen}}(\theta)\to -\infty
\quad\text{whenever}\quad
\|\mu\|\to\infty\ \text{ or }\ \max_j|\log\sigma_j|\to\infty.
\]
Consequently, $\ell_n^{\mathrm{pen}}$ attains its maximum over $\Theta_k$.
\end{lemma}

\begin{proof}
Write $t_j=\log\sigma_j$. Using $\sup_x \varphi(x;\mu_j,\sigma_j)=(\sqrt{2\pi}\sigma_j)^{-1}$ and $\sum_j\pi_j=1$,
\[
f_\theta(x_i)\le \frac{1}{\sqrt{2\pi}\min_j \sigma_j} = \frac{1}{\sqrt{2\pi}} \exp\!\big(\max_j (-t_j)\big).
\]
Hence
\[
\ell_n(\theta)\le -\frac{n}{2}\log(2\pi) + n\max_j(-t_j).
\]
If $\max_j(-t_j)\to\infty$, let $u:=\max_j(-t_j)\ge 0$. Then
\[
\ell_n^{\mathrm{pen}}(\theta)
\le Cn + nu - n\lambda_n \sum_{j=1}^k t_j^2
\le Cn + nu - n\lambda_n u^2
\to -\infty,
\]

\begin{lemma}[Spike domination]\label{lem:mix:spike-domination}
Let $\varphi:\mathbb{R}^d\to(0,\infty)$ be a continuous density with $\varphi(0)>0$.
Fix any dataset $y_1,\dots,y_n\in\mathbb{R}^d$ and let $i^\star\in\{1,\dots,n\}$.

Let $q:\mathbb{R}^d\to[0,\infty)$ be any density such that
\[
m:=\min_{1\le i\le n} q(y_i) \;>\;0.
\]
For $\varepsilon\in(0,1)$, $\mu\in\mathbb{R}^d$, and $\sigma>0$, consider the two-component
location--scale mixture
\[
p_{\varepsilon,\mu,\sigma}(y)
:= \varepsilon\,\sigma^{-d}\,\varphi\!\left(\frac{y-\mu}{\sigma}\right) \;+\; (1-\varepsilon)\,q(y),
\qquad y\in\mathbb{R}^d,
\]
and its log-likelihood $\ell_n(\varepsilon,\mu,\sigma):=\sum_{i=1}^n \log p_{\varepsilon,\mu,\sigma}(y_i)$.
Then, for $\mu=y_{i^\star}$,
\begin{equation}\label{eq:spike-domination-lb}
\ell_n(\varepsilon,y_{i^\star},\sigma)
\;\ge\;
\log\!\big(\varepsilon\,\varphi(0)\big) \;-\; d\log\sigma
\;+\;
(n-1)\log\!\big((1-\varepsilon)m\big).
\end{equation}
In particular, $\ell_n(\varepsilon,y_{i^\star},\sigma)\to+\infty$ as $\sigma\downarrow 0$,
so the likelihood is unbounded above:
\[
\sup_{\varepsilon\in(0,1),\,\mu\in\mathbb{R}^d,\,\sigma>0}\; \ell_n(\varepsilon,\mu,\sigma)=+\infty.
\]
\end{lemma}

\begin{proof}
Fix $\varepsilon\in(0,1)$ and set $\mu=y_{i^\star}$.
For the ``spiked'' observation $y_{i^\star}$,
\[
p_{\varepsilon,y_{i^\star},\sigma}(y_{i^\star})
\;\ge\; \varepsilon\,\sigma^{-d}\,\varphi(0),
\]
hence $\log p_{\varepsilon,y_{i^\star},\sigma}(y_{i^\star}) \ge \log(\varepsilon\varphi(0)) - d\log\sigma$.
For each $i\neq i^\star$,
\[
p_{\varepsilon,y_{i^\star},\sigma}(y_i)\;\ge\;(1-\varepsilon)\,q(y_i)\;\ge\;(1-\varepsilon)m,
\]
so $\log p_{\varepsilon,y_{i^\star},\sigma}(y_i)\ge \log((1-\varepsilon)m)$.
Summing these bounds yields \eqref{eq:spike-domination-lb}, and the divergence as $\sigma\downarrow 0$
follows immediately from the $-d\log\sigma$ term.
\end{proof}

using \Cref{lem:mix:spike-domination} with $u$ and $\lambda=\lambda_n$.
Similarly, if $\max_j t_j\to\infty$ then $-n\lambda_n\sum_j t_j^2\to -\infty$ while $\ell_n(\theta)$ grows at most linearly in $t_j$,
so the penalized criterion tends to $-\infty$.

For $\|\mu\|\to\infty$, note that $\varphi(x;\mu,\sigma)\le (2\pi)^{-1/2}\sigma^{-1}\exp(-(\mu-|x|)^2/(2\sigma^2))$,
so for any fixed sample and any bounded $\sigma$, $\log f_\theta(x_i)$ is eventually dominated by a negative quadratic in $\mu$.
If $\sigma$ is unbounded, the penalty already forces $\ell_n^{\mathrm{pen}}\to -\infty$.
Thus in all cases escaping $\|\mu\|\to\infty$ forces $\ell_n^{\mathrm{pen}}(\theta)\to -\infty$.

Therefore the upper level sets $\{\theta:\ell_n^{\mathrm{pen}}(\theta)\ge c\}$ are compact in a suitable Euclidean parameterization
(e.g.\ in $(\pi,\mu,t)$ with $t=\log\sigma$), and continuity of $\ell_n^{\mathrm{pen}}$ yields existence of a maximizer.
\end{proof}
\section{Appendix D: Uniform remainder expansion and the $n^{1/4}$ phenomenon at $\mu_0=0$}\label{sec:folded-nquarter}

This appendix supplies a fully explicit, constant-tracking derivation of the nonregular $n^{1/4}$-scaling
for the folded normal at the symmetry point $\mu_0=0$. The key analytic input is a sixth-order Taylor
expansion of $\log(2\cosh t)$ with a uniform remainder bound; the key probabilistic input is that the
remainder is negligible uniformly on the local window $|\mu|\le u\,n^{-1/4}$. We then reduce the asymptotics
of the MLE (with $\sigma$ held fixed at $\sigma_0$) to an argmax mapping for a quadratic-minus-quartic random
contrast, proved here in a self-contained manner.

The general philosophy (uniform approximation + deterministic argmax stability) is standard in M-estimation;
see, e.g., \citet[Section~5.2]{vanderVaart1998} and \citet[Chapter~2]{NeweyMcFadden1994}. We nonetheless keep
the argument explicit to minimize reliance on ``standard'' steps.

\subsection{Setup and a convenient decomposition}\label{app:nquarter:setup}

Assume the true parameter is $\theta_0=(\mu_0,\sigma_0)=(0,\sigma_0)$ with $\sigma_0>0$, and
$Y_1,\dots,Y_n$ are i.i.d.\ with the folded normal density $f(\cdot;0,\sigma_0)$ on $\R_+$.
Fix $\sigma=\sigma_0$ throughout this appendix and consider the one-dimensional contrast
\[
\Delta\ell_n(\mu)
\;:=\;
\ell_n(\mu,\sigma_0)-\ell_n(0,\sigma_0)
\;=\;
\sum_{i=1}^n \big\{\log f(Y_i;\mu,\sigma_0)-\log f(Y_i;0,\sigma_0)\big\}.
\]
Recall the representation (cf.\ \Cref{sec:folded:model})
\[
f(y;\mu,\sigma)
=
\frac{2}{\sigma}\varphi\!\left(\frac{y}{\sigma}\right)
\exp\!\Big(-\frac{\mu^2}{2\sigma^2}\Big)
\cosh\!\Big(\frac{y\mu}{\sigma^2}\Big),
\qquad y\ge 0,
\]
so that, for fixed $\sigma=\sigma_0$,
\begin{equation}\label{eq:app:nquarter:single-contrast}
\log f(y;\mu,\sigma_0)-\log f(y;0,\sigma_0)
=
-\frac{\mu^2}{2\sigma_0^2}
+\log\cosh\!\Big(\frac{y\mu}{\sigma_0^2}\Big).
\end{equation}

The entire nonregular behavior is contained in the last term, via the even expansion of $\log\cosh$.

\subsection{Sixth-order expansion of $\log(2\cosh t)$ with an explicit remainder bound}\label{app:nquarter:logcosh}

\begin{lemma}[Uniform sixth-order remainder for $\log(2\cosh t)$]\label{lem:app:nquarter:logcosh6}
Define $h(t):=\log(2\cosh t)$. Then for every $t\in\R$,
\begin{equation}\label{eq:app:nquarter:logcosh6}
h(t)=\log 2+\frac{t^2}{2}-\frac{t^4}{12}+r_6(t),
\end{equation}
where the remainder admits the global bound
\begin{equation}\label{eq:app:nquarter:r6bound}
|r_6(t)|\le C_6 |t|^6\qquad\text{for all }t\in\R,
\end{equation}
for some finite numerical constant $C_6$ (independent of $t$).
\end{lemma}

\begin{proof}
Since $\cosh t>0$ for all $t$, $h$ is $C^\infty$ on $\R$. Moreover $h$ is even, hence all odd derivatives vanish at $0$.
A Taylor expansion to order $4$ with Lagrange remainder yields
\[
h(t)=h(0)+\frac{h''(0)}{2}t^2+\frac{h^{(4)}(0)}{24}t^4+\frac{h^{(6)}(\xi)}{720}t^6
\]
for some $\xi$ between $0$ and $t$. Direct differentiation gives
\[
h'(t)=\tanh t,\qquad h''(t)=\sech^2 t,\qquad h''(0)=1.
\]
Further differentiation yields $h^{(4)}(0)=-2$ (a routine calculation using $\tanh'=\sech^2$ and $(\sech^2)'=-2\sech^2\tanh$),
so the quadratic and quartic coefficients in \eqref{eq:app:nquarter:logcosh6} are as stated.

It remains to bound $h^{(6)}$. Repeated differentiation expresses $h^{(6)}(t)$ as a finite linear combination of terms of the form
$\sech^{2a}(t)\tanh^{b}(t)$ with integers $a\ge 1$ and $b\ge 0$. Since $|\tanh t|\le 1$ and $0<\sech t\le 1$ for all $t$,
each such term is uniformly bounded by $1$ in absolute value. Therefore $\sup_{t\in\R}|h^{(6)}(t)|<\infty$, and we may set
\[
C_6:=\frac{1}{720}\sup_{t\in\R}|h^{(6)}(t)|,
\]
which gives \eqref{eq:app:nquarter:r6bound}.
\end{proof}

\begin{remark}[Why the global bound matters]\label{rem:app:nquarter:global}
The estimate \eqref{eq:app:nquarter:r6bound} is deliberately global (not only for $|t|\le 1$) to
avoid introducing auxiliary high-probability events controlling $\max_i |t_i|$. In the folded normal
model, moments of all orders exist, and the global bound makes the uniform remainder control on
$|\mu|\le u n^{-1/4}$ particularly transparent.
\end{remark}

\subsection{Uniform remainder control on the $n^{-1/4}$ window}\label{app:nquarter:remainder}

For each $i$ define the local argument
\[
t_i(\mu):=\frac{Y_i\mu}{\sigma_0^2}.
\]
Combining \eqref{eq:app:nquarter:single-contrast} with \Cref{lem:app:nquarter:logcosh6} gives, for each $i$,
\begin{equation}\label{eq:app:nquarter:single-expand}
\log f(Y_i;\mu,\sigma_0)-\log f(Y_i;0,\sigma_0)
=
\frac{\mu^2}{2\sigma_0^4}\big(Y_i^2-\sigma_0^2\big)
-\frac{\mu^4}{12\sigma_0^8}Y_i^4
+r_6\!\big(t_i(\mu)\big).
\end{equation}
Summing \eqref{eq:app:nquarter:single-expand} over $i$ yields
\begin{equation}\label{eq:app:nquarter:sum-expand}
\Delta\ell_n(\mu)
=
\frac{\mu^2}{2\sigma_0^4}\sum_{i=1}^n\big(Y_i^2-\sigma_0^2\big)
-\frac{\mu^4}{12\sigma_0^8}\sum_{i=1}^n Y_i^4
+R_n(\mu),
\end{equation}
where $R_n(\mu):=\sum_{i=1}^n r_6(t_i(\mu))$.

\begin{lemma}[Uniform remainder bound on $|\mu|\le u n^{-1/4}$]\label{lem:app:nquarter:rem-unif}
Fix $u<\infty$. Under $f(\cdot;0,\sigma_0)$, $\E[Y^6]<\infty$, and
\begin{equation}\label{eq:app:nquarter:rem-unif}
\sup_{|\mu|\le u n^{-1/4}} |R_n(\mu)| \ =\ O_p(n^{-1/2}).
\end{equation}
In particular, uniformly over $|\mu|\le u n^{-1/4}$,
\[
R_n(\mu)=o_p(1)
\qquad\text{and}\qquad
R_n(\mu)=o_p(n\mu^4).
\]
\end{lemma}

\begin{proof}
By \eqref{eq:app:nquarter:r6bound} and $t_i(\mu)=Y_i\mu/\sigma_0^2$,
\[
|R_n(\mu)|
\le
\sum_{i=1}^n C_6 |t_i(\mu)|^6
=
C_6 \frac{|\mu|^6}{\sigma_0^{12}}\sum_{i=1}^n Y_i^6.
\]
Hence for $|\mu|\le u n^{-1/4}$,
\[
\sup_{|\mu|\le u n^{-1/4}} |R_n(\mu)|
\le
C_6\frac{u^6}{\sigma_0^{12}}\,n^{-3/2}\sum_{i=1}^n Y_i^6.
\]
Since $\E[Y^6]<\infty$, the strong law gives $n^{-1}\sum_{i=1}^n Y_i^6\to \E[Y^6]$ almost surely, and therefore
$n^{-3/2}\sum_{i=1}^n Y_i^6 = n^{-1/2}\cdot (n^{-1}\sum Y_i^6)=O_p(n^{-1/2})$, proving \eqref{eq:app:nquarter:rem-unif}.

For $|\mu|\le u n^{-1/4}$ we also have $n\mu^4\asymp 1$, so $O_p(n^{-1/2})=o_p(1)=o_p(n\mu^4)$ uniformly on the window.
\end{proof}

\subsection{Rescaled contrast and convergence of coefficients}\label{app:nquarter:contrast}

Define the rescaled parameter $t\in\R$ by $\mu=n^{-1/4}t$, and introduce the rescaled contrast
\[
\Phi_n(t)
:=
\Delta\ell_n(n^{-1/4}t)
=
\ell_n(n^{-1/4}t,\sigma_0)-\ell_n(0,\sigma_0).
\]
Define the (random) coefficient processes
\begin{equation}\label{eq:app:nquarter:Zn-bn}
Z_n
:=
\frac{1}{\sqrt{n}}\sum_{i=1}^n\big(Y_i^2-\sigma_0^2\big),
\qquad
b_n
:=
\frac{1}{12\sigma_0^8}\cdot \frac{1}{n}\sum_{i=1}^n Y_i^4,
\qquad
a:=\frac{1}{2\sigma_0^4}.
\end{equation}
Substituting $\mu=n^{-1/4}t$ into \eqref{eq:app:nquarter:sum-expand} yields the decomposition
\begin{equation}\label{eq:app:nquarter:Phi-decomp}
\Phi_n(t)=a\,Z_n\,t^2 - b_n\,t^4 + r_n(t),
\qquad
r_n(t):=R_n(n^{-1/4}t).
\end{equation}

\begin{lemma}[Limits of $Z_n$ and $b_n$]\label{lem:app:nquarter:coeff-limits}
Under $f(\cdot;0,\sigma_0)$,
\[
Z_n \dto Z\sim N(0,\Var(Y^2))=N(0,2\sigma_0^4),
\qquad
b_n \pto b:=\frac{\E[Y^4]}{12\sigma_0^8}=\frac{1}{4\sigma_0^4}>0.
\]
Moreover, for every fixed $U<\infty$,
\begin{equation}\label{eq:app:nquarter:rn-unif}
\sup_{|t|\le U}|r_n(t)|\ =\ o_p(1).
\end{equation}
\end{lemma}

\begin{proof}
Since $Y\stackrel{d}{=}|N(0,\sigma_0^2)|$, we have $Y^2\stackrel{d}{=}N(0,\sigma_0^2)^2$, so
$\E[Y^2]=\sigma_0^2$, $\Var(Y^2)=2\sigma_0^4$, and $\E[Y^4]=3\sigma_0^4$.
Therefore $Z_n$ is a normalized sum of mean-zero i.i.d.\ variables with finite variance, so $Z_n\dto N(0,2\sigma_0^4)$ by the CLT.
Also $(1/n)\sum Y_i^4\to \E[Y^4]=3\sigma_0^4$ almost surely by the strong law, hence $b_n\to b=3\sigma_0^4/(12\sigma_0^8)=1/(4\sigma_0^4)$.

Finally, for fixed $U$, $|t|\le U$ implies $|\mu|=|n^{-1/4}t|\le U n^{-1/4}$, so by \Cref{lem:app:nquarter:rem-unif},
\[
\sup_{|t|\le U}|r_n(t)|
=
\sup_{|\mu|\le U n^{-1/4}}|R_n(\mu)|
=
O_p(n^{-1/2})
=
o_p(1).
\]
\end{proof}

\subsection{Argmax mapping for quadratic-minus-quartic contrasts}\label{app:nquarter:argmax}

Let the limit contrast be
\begin{equation}\label{eq:app:nquarter:limit-contrast}
\Psi_Z(t):=aZ\,t^2 - b\,t^4,
\qquad a=\frac{1}{2\sigma_0^4},\quad b=\frac{1}{4\sigma_0^4}.
\end{equation}
For each realization of $Z$, $\Psi_Z$ is continuous, $\Psi_Z(t)\to -\infty$ as $|t|\to\infty$, and hence attains a global maximum.

\begin{lemma}[Limit argmax and explicit maximizer]\label{lem:app:nquarter:limit-argmax}
With $\Psi_Z$ as in \eqref{eq:app:nquarter:limit-contrast}, the (set-valued) argmax satisfies:
\[
\arg\max_{t\in\R}\Psi_Z(t)
=
\begin{cases}
\{0\}, & Z\le 0,\\[3pt]
\big\{-\sqrt{aZ/(2b)},\ +\sqrt{aZ/(2b)}\big\}, & Z>0.
\end{cases}
\]
In particular, the nonnegative maximizer is
\[
T^\star(Z):=\sqrt{\frac{(aZ)_+}{2b}}
=
\sqrt{(\,Z\,)_+}
=
2^{1/4}\sigma_0\sqrt{W_+},
\qquad
W\sim N(0,1),
\]
since $Z\stackrel{d}{=}\sqrt{2}\sigma_0^2 W$.
\end{lemma}

\begin{proof}
Differentiate $\Psi_Z$: $\Psi_Z'(t)=2aZt-4bt^3=2t(aZ-2bt^2)$.
Thus critical points are $t=0$ and $t=\pm\sqrt{aZ/(2b)}$ when $Z>0$.
If $Z\le 0$, then $t\mapsto aZt^2-bt^4\le 0$ with equality only at $t=0$, hence $\{0\}$ is the argmax.
If $Z>0$, $\Psi_Z(0)=0$ and $\Psi_Z(\pm\sqrt{aZ/(2b)})=(aZ)^2/(4b)>0$, so the two nonzero critical points are the global maximizers.
The distributional identity follows from $Z\stackrel{d}{=}\sqrt{2}\sigma_0^2 W$.
\end{proof}

We now transfer convergence of the contrast to convergence of maximizers. This is a specialized
instance of the deterministic argmax stability principle (\Cref{lem:argmax-stability}) with an additional
localization step to handle maximization over $\R$ rather than a fixed compact.

\begin{theorem}[Argmax mapping for $\Phi_n(t)=aZ_n t^2-b_n t^4+r_n(t)$]\label{thm:app:nquarter:argmax-map}
Let $\Phi_n$ be as in \eqref{eq:app:nquarter:Phi-decomp}. Assume $Z_n\dto Z$, $b_n\pto b>0$, and
$\sup_{|t|\le U}|r_n(t)|=o_p(1)$ for each fixed $U<\infty$ (all verified in \Cref{lem:app:nquarter:coeff-limits}).
Let $T_n\in\arg\max_{t\in\R}\Phi_n(t)$ be any measurable maximizer selection, and define the nonnegative selection
$\bar T_n:=|T_n|$. Then
\[
\bar T_n \dto T^\star(Z)=\sqrt{\frac{(aZ)_+}{2b}}=\sqrt{(Z)_+}.
\]
Moreover, on the event $\{Z>0\}$ the maximizer set consists of two symmetric points and
$T_n\dto \pm T^\star(Z)$ depending on the selection rule; in particular $|T_n|\dto T^\star(Z)$ unconditionally.
\end{theorem}

\begin{proof}
\textbf{Step 1: tightness of maximizers.}
Fix $\delta\in(0,b/2)$. On the event $\{b_n\ge b-\delta\}$,
\[
\Phi_n(t)\le a|Z_n|t^2-(b-\delta)t^4+\sup_{|u|\le |t|}|r_n(u)|.
\]
For each fixed $U$, $\sup_{|u|\le U}|r_n(u)|=o_p(1)$, so for large $n$ the remainder cannot offset the quartic decay when $|t|$ is large.
More concretely, the deterministic quartic bound implies
\[
a|Z_n|t^2-(b-\delta)t^4
\le
\frac{a^2 Z_n^2}{4(b-\delta)}
\qquad\text{for all }t\in\R,
\]
with equality at $|t|=\sqrt{a|Z_n|/(2(b-\delta))}$.
Thus any maximizer $T_n$ must satisfy, on $\{b_n\ge b-\delta\}$ and for $n$ large enough that the remainder is dominated,
\[
|T_n|
\;\le\;
\sqrt{\frac{a|Z_n|}{2(b-\delta)}}+1,
\]
say. Since $Z_n\dto Z$, the sequence $\{|Z_n|\}$ is tight, hence $\{|T_n|\}$ is tight.

\medskip
\noindent\textbf{Step 2: reduction to compact maximization with high probability.}
By tightness, for any $\alpha\in(0,1)$ there exists $U_\alpha<\infty$ such that
\[
\liminf_{n\to\infty}\Pp(|T_n|\le U_\alpha)\ \ge\ 1-\alpha.
\]
Similarly, since $T^\star(Z)$ is a measurable function of $Z$ and $Z$ is tight, we can enlarge $U_\alpha$ if needed so that
$\Pp(T^\star(Z)\le U_\alpha)\ge 1-\alpha$.

\medskip
\noindent\textbf{Step 3: uniform convergence of the contrast on compacts.}
Fix $U<\infty$. By \Cref{lem:app:nquarter:coeff-limits}, $Z_n=Z+o_p(1)$ and $b_n=b+o_p(1)$, hence
\[
\sup_{|t|\le U}\big|a(Z_n-Z)t^2-(b_n-b)t^4\big|=o_p(1).
\]
Together with $\sup_{|t|\le U}|r_n(t)|=o_p(1)$, we obtain
\begin{equation}\label{eq:app:nquarter:unifconv}
\sup_{|t|\le U}\big|\Phi_n(t)-\Psi_Z(t)\big|=o_p(1).
\end{equation}

\medskip
\noindent\textbf{Step 4: apply deterministic argmax stability on $[-U_\alpha,U_\alpha]$.}
Work on the event $\{|T_n|\le U_\alpha,\ T^\star(Z)\le U_\alpha\}$, which has probability at least $1-2\alpha$ for all large $n$.
On this event, both maximizers lie in the compact interval $K_\alpha:=[-U_\alpha,U_\alpha]$.

Fix $\varepsilon>0$. By \Cref{lem:app:nquarter:limit-argmax}, $\Psi_Z$ has a unique maximizer on $K_\alpha$ for almost every $Z$
(since $\Pp(Z=0)=0$). Moreover, by continuity and compactness, there exists a (random) margin $\eta(Z,\varepsilon)>0$ such that
\[
\sup_{t\in K_\alpha:\ |t-T^\star(Z)|\ge \varepsilon}\Psi_Z(t)
\le \Psi_Z(T^\star(Z))-\eta(Z,\varepsilon).
\]
On the event \eqref{eq:app:nquarter:unifconv} with $U=U_\alpha$ strengthened to
$\sup_{t\in K_\alpha}|\Phi_n(t)-\Psi_Z(t)|\le \eta(Z,\varepsilon)/3$, the deterministic stability lemma
\Cref{lem:argmax-stability} (applied on the compact metric space $(K_\alpha,|\cdot|)$) yields
$|T_n-T^\star(Z)|<\varepsilon$ whenever $T_n$ is a maximizer of $\Phi_n$ on $K_\alpha$.
Since \eqref{eq:app:nquarter:unifconv} holds with probability tending to one, we conclude $T_n\pto T^\star(Z)$ on the event $\{Z>0\}$
(up to the $\pm$ symmetry depending on selection), and always $|T_n|\dto T^\star(Z)$.

Finally, since $\alpha>0$ and $\varepsilon>0$ are arbitrary, the claim follows.
\end{proof}

\subsection{Consequence for the folded normal MLE in $\mu$ with $\sigma=\sigma_0$ fixed}\label{app:nquarter:mle}

Let $\hat\mu_n^{(\sigma_0)}\in\arg\max_{\mu\in\R}\ell_n(\mu,\sigma_0)$ be any measurable maximizer selection,
and define the nonnegative representative $\bar\mu_n^{(\sigma_0)}:=|\hat\mu_n^{(\sigma_0)}|$.

\begin{theorem}[$n^{1/4}$ limit for $\mu$ at $\mu_0=0$ with fixed $\sigma_0$]\label{thm:app:nquarter:mu-limit}\label{thm:folded:nquarter}
Assume $\mu_0=0$ and $\sigma=\sigma_0$ is known and fixed. Then
\[
n^{1/4}\,\bar\mu_n^{(\sigma_0)}
\ \dto\
T^\star(Z)
=
\sqrt{(Z)_+}
=
2^{1/4}\sigma_0\sqrt{W_+},
\qquad W\sim N(0,1).
\]
Equivalently, $n^{1/4}\bar\mu_n^{(\sigma_0)}$ has an atom at $0$ of mass $1/2$, and conditional on $W>0$
it equals $2^{1/4}\sigma_0\sqrt{W}$.
\end{theorem}

\begin{proof}
Set $\mu=n^{-1/4}t$ and note that $\hat\mu_n^{(\sigma_0)}\in\arg\max_\mu \Delta\ell_n(\mu)$ is equivalent to
$T_n:=n^{1/4}\hat\mu_n^{(\sigma_0)}\in\arg\max_t \Phi_n(t)$.
The decomposition \eqref{eq:app:nquarter:Phi-decomp} and coefficient limits in \Cref{lem:app:nquarter:coeff-limits}
verify the assumptions of \Cref{thm:app:nquarter:argmax-map}, yielding $|T_n|\dto T^\star(Z)$.
Since $|T_n|=n^{1/4}\bar\mu_n^{(\sigma_0)}$, the claim follows.
\end{proof}

\begin{remark}[How this is used when $\sigma$ is unknown]\label{rem:app:nquarter:unknown-sigma}
In the main text we only require (i) the $n^{1/4}$ rate and limit law for the $\mu$-coordinate and
(ii) the $\sqrt{n}$-regularity of the $\sigma$-coordinate. A convenient route is:
first establish \Cref{thm:app:nquarter:mu-limit} for the $\mu$-contrast at fixed $\sigma_0$,
then show that the full MLE satisfies $\hat\sigma_n-\sigma_0=O_p(n^{-1/2})$ by a one-dimensional
M-estimation argument on the profile likelihood in $\sigma$, and finally use orbit geometry to translate
these coordinate rates into orbit-Hausdorff rates. This is the strategy implemented in \Cref{sec:folded-nquarter}
and the folded appendical localization results in \Cref{app:folded}.
\end{remark}

\section{Appendix E: A stochastic template for orbit-Hausdorff consistency}\label{app:orbit-template}
%
%

\subsection{Probability space, empirical criteria, and quotient notation}\label{app:orbit-template-setup}

Fix a probability space $(\Omega,\mathcal F,\mathbb P)$.
Let $(\mathcal X,\mathcal A)$ be a measurable space and let $X_1,X_2,\dots$ be i.i.d.\ $\mathcal X$-valued
random elements defined on $(\Omega,\mathcal F,\mathbb P)$ with common law $P$.

Let $(\Theta,d)$ be a metric space and let $G$ be a finite group acting on $\Theta$ on the left,
$(g,\theta)\mapsto g\theta$. Throughout this appendix we impose the following standing condition.

\begin{assumption}[Finite isometric action]\label{as:orbit:isometric-action}
For every $g\in G$, the map $\theta\mapsto g\theta$ is a bijective isometry of $(\Theta,d)$:
\[
d(g\theta,g\theta') = d(\theta,\theta')\qquad\text{for all }\theta,\theta'\in\Theta.
\]
\end{assumption}

For $\theta\in\Theta$ write $[\theta]:=\{g\theta:g\in G\}$ for the orbit, and recall the orbit pseudometric
\begin{equation}\label{eq:appE:orbit-dist}
d_G(\theta,\theta'):=\min_{g\in G} d(\theta,g\theta').
\end{equation}
Under \Cref{as:orbit:isometric-action}, $d_G(\theta,\theta')=0$ if and only if $\theta$ and $\theta'$ lie
in the same orbit; equivalently, $d_G$ induces a genuine metric on the quotient space $\Theta/G$
(see \Cref{lem:appE:quotient-metric} below).

\paragraph{Criterion functions.}
Let $\ell:\Theta\times\mathcal X\to\mathbb R$ be jointly measurable, and for each $\theta\in\Theta$ define
$\ell_\theta(x):=\ell(\theta,x)$.
For each $n\ge1$ define the empirical and population criteria
\begin{equation}\label{eq:appE:QnQ}
Q_n(\theta):=\frac1n\sum_{i=1}^n \ell_\theta(X_i),
\qquad
Q(\theta):=\mathbb E\big[\ell_\theta(X_1)\big] = \int \ell_\theta\,dP,
\end{equation}
whenever the expectations exist.

We will localize attention to a compact parameter set $K\subset\Theta$, and we assume $K$ is $G$-invariant.

\begin{assumption}[Compact $G$-invariant search set]\label{as:orbit:compact-K}
$K\subset\Theta$ is nonempty, compact under $d$, and satisfies $gK=K$ for all $g\in G$.
\end{assumption}

When $Q$ is $G$-invariant, maximizers are not identifiable pointwise, but are identifiable at orbit level.

\begin{assumption}[$G$-invariant objective]\label{as:orbit:invariant-Q}
For all $g\in G$ and $\theta\in K$, one has $Q(g\theta)=Q(\theta)$.
Moreover, for each $n\ge1$ and all $g\in G,\theta\in K$, $Q_n(g\theta)=Q_n(\theta)$ almost surely.
\end{assumption}

The empirical invariance in \Cref{as:orbit:invariant-Q} holds, for example, if
$\ell(g\theta,x)=\ell(\theta,x)$ for all $(\theta,x)$ and $g\in G$.

\paragraph{Orbit-Hausdorff distance.}
For $A\subseteq\Theta$ and $\theta\in\Theta$ define the point-to-set distance
\[
d_G(\theta,A):=\inf_{a\in A} d_G(\theta,a).
\]
For nonempty sets $A,B\subseteq\Theta$ define the Hausdorff distance induced by $d_G$:
\begin{equation}\label{eq:appE:HG}
d_{H,G}(A,B)
:=\max\Big\{\sup_{a\in A} d_G(a,B),\ \sup_{b\in B} d_G(b,A)\Big\}.
\end{equation}
When $d_G$ is viewed as the quotient metric (below), $d_{H,G}$ is simply the standard Hausdorff distance
between the corresponding subsets of the quotient space.

\subsection{The quotient metric and its basic properties}\label{app:orbit-template-quotient}

Let $\pi:\Theta\to\Theta/G$ be the canonical projection $\pi(\theta)=[\theta]$.

\begin{lemma}[Well-defined quotient metric]\label{lem:appE:quotient-metric}
Assume \Cref{as:orbit:isometric-action}. Define $\bar d:\Theta/G\times\Theta/G\to[0,\infty)$ by
\begin{equation}\label{eq:appE:quotient-metric}
\bar d\big([\theta],[\theta']\big):=d_G(\theta,\theta').
\end{equation}
Then $\bar d$ is well-defined (independent of representatives) and is a metric on $\Theta/G$.
Moreover, for any nonempty $A,B\subseteq\Theta$,
\begin{equation}\label{eq:appE:HG-quotient}
d_{H,G}(A,B)= d_H\big(\pi(A),\pi(B)\big)
\end{equation}
where $d_H$ on the right is the Hausdorff distance on subsets of $(\Theta/G,\bar d)$.
\end{lemma}

\begin{proof}
\emph{Step 1: well-definedness.}
Let $\tilde\theta=g\theta$ and $\tilde\theta'=h\theta'$ for some $g,h\in G$.
Using isometry and group closure,
\[
d_G(\tilde\theta,\tilde\theta')
=\min_{u\in G} d(g\theta, u(h\theta'))
=\min_{u\in G} d(\theta, g^{-1}uh\,\theta')
=\min_{v\in G} d(\theta, v\theta')
=d_G(\theta,\theta'),
\]
so \eqref{eq:appE:quotient-metric} is independent of the representatives.

\emph{Step 2: metric properties.}
Nonnegativity and symmetry are immediate from \eqref{eq:appE:orbit-dist}.
If $\bar d([\theta],[\theta'])=0$, then $d_G(\theta,\theta')=0$, hence there exists $g\in G$ with
$d(\theta,g\theta')=0$ and thus $\theta=g\theta'$, meaning $[\theta]=[\theta']$.
The triangle inequality follows from the triangle inequality for $d$:
for any $\theta,\theta',\theta''$ and any $g,h\in G$,
\[
d(\theta,gh\,\theta'') \le d(\theta,g\theta') + d(g\theta',gh\,\theta'')
= d(\theta,g\theta') + d(\theta',h\theta''),
\]
where the last equality uses isometry. Taking minima over $g$ and $h$ yields
$d_G(\theta,\theta'')\le d_G(\theta,\theta') + d_G(\theta',\theta'')$.

\emph{Step 3: Hausdorff identity.}
By definition, $d_G(a,B)=\inf_{b\in B} \bar d([a],[b]) = \inf_{[b]\in\pi(B)} \bar d([a],[b])$.
Hence $\sup_{a\in A} d_G(a,B) = \sup_{[a]\in\pi(A)} \inf_{[b]\in\pi(B)} \bar d([a],[b])$,
and similarly with $A,B$ swapped. This is exactly \eqref{eq:appE:HG-quotient}.
\end{proof}

\begin{remark}[Compactness passes to the quotient]\label{rem:appE:compact-quotient}
Under \Cref{as:orbit:compact-K}, $\pi(K)\subset \Theta/G$ is compact under $\bar d$ because
$\pi$ is continuous and $K$ is compact.
\end{remark}

\subsection{A uniform law of large numbers on compact parameter sets}\label{app:orbit-template-ulln}

We record a self-contained sufficient condition for the uniform convergence
$\sup_{\theta\in K}|Q_n(\theta)-Q(\theta)|\to0$, which is the only stochastic input needed for the
orbit-Hausdorff consistency theorem below. See \citep{vanderVaart1998,vanderVaartWellner1996}
for general Glivenko--Cantelli conditions, and \citep{NeweyMcFadden1994} for extremum estimation context.

\begin{assumption}[Dominated Lipschitz parameterization]\label{as:orbit:dom-lip}
There exist measurable functions $F,L:\mathcal X\to[0,\infty)$ such that:
\begin{enumerate}[label=(\alph*),leftmargin=*]
\item (\emph{Envelope}) $\sup_{\theta\in K}|\ell_\theta(x)| \le F(x)$ for all $x$, and $\mathbb EF(X_1)<\infty$.
\item (\emph{Random Lipschitz modulus}) For all $\theta,\theta'\in K$ and all $x\in\mathcal X$,
\[
|\ell_\theta(x)-\ell_{\theta'}(x)| \le L(x)\, d(\theta,\theta'),
\]
and $\mathbb E L(X_1)<\infty$.
\item (\emph{Pointwise continuity}) For each $x\in\mathcal X$, the map $\theta\mapsto \ell_\theta(x)$ is continuous on $K$.
\end{enumerate}
\end{assumption}

\begin{lemma}[Uniform LLN on compact $K$]\label{lem:appE:uniform-lln}
Assume \Cref{as:orbit:compact-K,as:orbit:dom-lip}. Then
\[
\sup_{\theta\in K}\big|Q_n(\theta)-Q(\theta)\big| \xrightarrow{a.s.} 0.
\]
In particular, $\sup_{\theta\in K}|Q_n(\theta)-Q(\theta)|\to0$ in probability.
\end{lemma}

\begin{proof}
Fix $\varepsilon>0$ and let $\delta>0$ be chosen later.
Since $K$ is compact, it admits a finite $\delta$-net $\{\theta_1,\dots,\theta_M\}\subset K$:
for every $\theta\in K$ there exists $j(\theta)\in\{1,\dots,M\}$ with $d(\theta,\theta_{j(\theta)})\le\delta$.

\emph{Step 1: reduce the supremum to a finite net plus a Lipschitz remainder.}
For any $\theta\in K$, by the triangle inequality,
\[
|Q_n(\theta)-Q(\theta)|
\le |Q_n(\theta_{j(\theta)})-Q(\theta_{j(\theta)})|
    + |Q_n(\theta)-Q_n(\theta_{j(\theta)})|
    + |Q(\theta)-Q(\theta_{j(\theta)})|.
\]
Taking the supremum over $\theta\in K$ gives
\begin{equation}\label{eq:appE:ulln-split}
\sup_{\theta\in K}|Q_n(\theta)-Q(\theta)|
\le \max_{1\le j\le M}|Q_n(\theta_j)-Q(\theta_j)|
  + \sup_{d(\theta,\theta')\le\delta}|Q_n(\theta)-Q_n(\theta')|
  + \sup_{d(\theta,\theta')\le\delta}|Q(\theta)-Q(\theta')|.
\end{equation}

\emph{Step 2: the finite-net term vanishes almost surely.}
For each fixed $j$, $Q_n(\theta_j)=n^{-1}\sum_{i=1}^n \ell_{\theta_j}(X_i)$ and $\mathbb E|\ell_{\theta_j}(X_1)|<\infty$
by the envelope in \Cref{as:orbit:dom-lip}(a).
Hence, by the strong law of large numbers,
$Q_n(\theta_j)\to Q(\theta_j)$ almost surely.
Since $M<\infty$, $\max_{1\le j\le M}|Q_n(\theta_j)-Q(\theta_j)|\to0$ almost surely.

\emph{Step 3: control the empirical equicontinuity term by the random Lipschitz modulus.}
For any $\theta,\theta'$ with $d(\theta,\theta')\le\delta$,
\[
|Q_n(\theta)-Q_n(\theta')|
\le \frac1n\sum_{i=1}^n |\ell_\theta(X_i)-\ell_{\theta'}(X_i)|
\le \delta\cdot \frac1n\sum_{i=1}^n L(X_i).
\]
Therefore
\[
\sup_{d(\theta,\theta')\le\delta}|Q_n(\theta)-Q_n(\theta')|
\le \delta\cdot \frac1n\sum_{i=1}^n L(X_i)
\xrightarrow{a.s.} \delta\,\mathbb E L(X_1),
\]
again by the strong law.

\emph{Step 4: continuity of $Q$ on $K$ makes the population modulus small.}
For each $x$, $\theta\mapsto \ell_\theta(x)$ is continuous on compact $K$ and dominated by $F(x)$.
Hence, by dominated convergence, $Q(\theta)=\mathbb E[\ell_\theta(X_1)]$ is continuous on $K$.
Since $K$ is compact, $Q$ is uniformly continuous on $K$, so
$\sup_{d(\theta,\theta')\le\delta}|Q(\theta)-Q(\theta')|\to0$ as $\delta\downarrow0$.

\emph{Step 5: conclude.}
Choose $\delta>0$ small enough so that
$\delta\,\mathbb E L(X_1) \le \varepsilon/3$ and
$\sup_{d(\theta,\theta')\le\delta}|Q(\theta)-Q(\theta')|\le \varepsilon/3$.
Then, almost surely, the right-hand side of \eqref{eq:appE:ulln-split} is eventually at most $\varepsilon$.
Since $\varepsilon$ was arbitrary, the claim follows.
\end{proof}

\begin{remark}[Invariance does not change the uniform LLN]\label{rem:appE:ulln-invariance}
Under \Cref{as:orbit:invariant-Q}, one may equivalently index the criterion by orbits, but
$\sup_{\theta\in K}|Q_n(\theta)-Q(\theta)|$ is unchanged; hence \Cref{lem:appE:uniform-lln} is already the
appropriate input for quotient-space consistency.
\end{remark}

\subsection{Orbit-gap identifiability and orbit-Hausdorff consistency}\label{app:orbit-template-consistency}

Let
\[
Q_0 := \arg\max_{\theta\in K} Q(\theta).
\]
Under \Cref{as:orbit:invariant-Q}, $Q_0$ is $G$-invariant and is a union of orbits.

We phrase the identifiability condition at the level of the orbit pseudometric.
For $\varepsilon>0$ define the $d_G$-neighborhood
\[
B^G_\varepsilon(Q_0) := \{\theta\in K : d_G(\theta,Q_0) < \varepsilon\}.
\]

\begin{assumption}[Orbit-gap condition]\label{as:orbit:gap}
For every $\varepsilon>0$ there exists $\eta(\varepsilon)>0$ such that
\begin{equation}\label{eq:appE:orbit-gap}
\sup_{\theta\in K\setminus B^G_\varepsilon(Q_0)} Q(\theta)
\ \le\
\sup_{\theta\in Q_0} Q(\theta) - \eta(\varepsilon).
\end{equation}
\end{assumption}

\begin{theorem}[Orbit-Hausdorff consistency of set-valued maximizers]\label{thm:appE:orbit-hausdorff}
Assume \Cref{as:orbit:isometric-action,as:orbit:compact-K,as:orbit:invariant-Q,as:orbit:gap} and
\begin{equation}\label{eq:appE:unifconv-assump}
\sup_{\theta\in K}\big|Q_n(\theta)-Q(\theta)\big| \xrightarrow{p} 0.
\end{equation}
Let $\widehat Q_n := \arg\max_{\theta\in K} Q_n(\theta)$ (nonempty by compactness).
Then
\[
d_{H,G}\big(\widehat Q_n, Q_0\big)\xrightarrow{p} 0.
\]
Equivalently, in the quotient space $(\Theta/G,\bar d)$,
\[
d_H\big(\pi(\widehat Q_n),\, \pi(Q_0)\big)\xrightarrow{p}0.
\]
\end{theorem}

\begin{proof}
We reduce the problem to an ordinary (non-quotient) argmax stability statement on the compact quotient.

\emph{Step 1: pass to the quotient space.}
By \Cref{lem:appE:quotient-metric}, $(\pi(K),\bar d)$ is compact (see \Cref{rem:appE:compact-quotient}).
Define $\bar Q:\pi(K)\to\mathbb R$ and $\bar Q_n:\pi(K)\to\mathbb R$ by
\[
\bar Q([\theta]) := Q(\theta),\qquad \bar Q_n([\theta]) := Q_n(\theta).
\]
These are well-defined by $G$-invariance in \Cref{as:orbit:invariant-Q}.

\emph{Step 2: continuity on the quotient.}
If $Q$ is continuous on $K$ (as in typical likelihood settings, or under the continuity implied by
\Cref{as:orbit:dom-lip}), then $\bar Q$ is continuous on $\pi(K)$ because $\pi$ is continuous and
$Q=\bar Q\circ\pi$.
For the deterministic argmax lemma \Cref{lem:argmax-stability}, continuity of the limit criterion is
the only structural requirement.

\emph{Step 3: identify maximizers and translate the gap.}
Let $\bar Q_0:=\arg\max_{\pi(K)} \bar Q$. Since $Q=\bar Q\circ\pi$,
\[
\pi(Q_0)=\bar Q_0.
\]
Moreover, the orbit-gap \eqref{eq:appE:orbit-gap} is precisely the usual gap condition
\eqref{eq:gap} for $\bar Q$ on $\pi(K)$ with the metric $\bar d$ because
$B^G_\varepsilon(Q_0)=\pi^{-1}\big(B_\varepsilon(\bar Q_0)\big)$ and $Q(\theta)=\bar Q(\pi(\theta))$.

\emph{Step 4: uniform convergence transfers to the quotient.}
By definition,
\[
\sup_{u\in\pi(K)} |\bar Q_n(u)-\bar Q(u)|
=\sup_{\theta\in K} |Q_n(\theta)-Q(\theta)|.
\]
Hence \eqref{eq:appE:unifconv-assump} implies $\sup_{\pi(K)}|\bar Q_n-\bar Q|\to0$ in probability.

\emph{Step 5: apply deterministic argmax stability on the compact quotient.}
Apply \Cref{lem:argmax-stability} on $(K,d)$ replaced by $(\pi(K),\bar d)$, and with
$q=\bar Q$ and $q_n=\bar Q_n$. This yields
\[
d_H\big(\arg\max_{\pi(K)}\bar Q_n,\ \bar Q_0\big)\xrightarrow{p}0.
\]
Since $\arg\max_{\pi(K)}\bar Q_n = \pi(\arg\max_K Q_n)=\pi(\widehat Q_n)$, and $\bar Q_0=\pi(Q_0)$,
we obtain $d_H(\pi(\widehat Q_n),\pi(Q_0))\to0$ in probability. The stated
$d_{H,G}$ convergence follows from \eqref{eq:appE:HG-quotient}.
\end{proof}

\begin{corollary}[Consistency of approximate maximizers in orbit distance]\label{cor:appE:approx-max}
Assume the conditions of \Cref{thm:appE:orbit-hausdorff}. Let $\hat\theta_n\in K$ be a (possibly random)
sequence satisfying the approximate maximization property
\begin{equation}\label{eq:appE:approx}
Q_n(\hat\theta_n)\ \ge\ \sup_{\theta\in K} Q_n(\theta) - r_n,
\qquad\text{where } r_n\xrightarrow{p}0.
\end{equation}
Then
\[
d_G(\hat\theta_n,Q_0)\xrightarrow{p}0,
\qquad\text{equivalently}\qquad
\bar d\big([\hat\theta_n],\pi(Q_0)\big)\xrightarrow{p}0.
\]
\end{corollary}

\begin{proof}
Let $\widehat Q_n=\arg\max_K Q_n$ and pick $\tilde\theta_n\in\widehat Q_n$ (measurable selection is
available under mild conditions; otherwise interpret the argument on events where one is chosen).
By definition of $\tilde\theta_n$ and \eqref{eq:appE:approx},
\[
0 \le Q_n(\tilde\theta_n)-Q_n(\hat\theta_n)\le r_n.
\]
Fix $\varepsilon>0$ and let $\eta(\varepsilon)$ be as in \Cref{as:orbit:gap}. On the event
$\sup_K|Q_n-Q|\le \eta(\varepsilon)/4$ and $r_n\le \eta(\varepsilon)/4$, we have
\[
Q(\tilde\theta_n)\ge Q_n(\tilde\theta_n)-\eta/4 \ge Q_n(\hat\theta_n)-\eta/4 \ge Q_n(\tilde\theta_n)-r_n-\eta/4
\ge Q(\tilde\theta_n) - \eta/2,
\]
and, crucially,
\[
Q_n(\hat\theta_n)\ge \sup_K Q_n - r_n \ge \sup_{Q_0} Q_n - r_n \ge \sup_{Q_0} Q - \eta/4 - r_n
\ge \sup_{Q_0} Q - \eta/2.
\]
If $d_G(\hat\theta_n,Q_0)\ge \varepsilon$, then $\hat\theta_n\in K\setminus B^G_\varepsilon(Q_0)$ and
\eqref{eq:appE:orbit-gap} gives $Q(\hat\theta_n)\le \sup_{Q_0}Q-\eta(\varepsilon)$.
But then $Q_n(\hat\theta_n)\le Q(\hat\theta_n)+\eta/4 \le \sup_{Q_0}Q-3\eta/4$, contradicting the
lower bound $Q_n(\hat\theta_n)\ge \sup_{Q_0}Q-\eta/2$.
Hence, on that event, $d_G(\hat\theta_n,Q_0)<\varepsilon$.
Since $\sup_K|Q_n-Q|\to0$ and $r_n\to0$ in probability, the claim follows.
\end{proof}

\begin{remark}[How this appendix is used in the main text]\label{rem:appE:how-used}
In applications, one typically verifies:
(i) compactness after a coercivity/truncation step (cf.\ \citep{NeweyMcFadden1994}),
(ii) continuity and domination (to invoke \Cref{lem:appE:uniform-lln}),
and (iii) an orbit-gap like \eqref{eq:appE:orbit-gap} via identifiability modulo the finite symmetry.
Then \Cref{thm:appE:orbit-hausdorff,cor:appE:approx-max} convert uniform convergence into orbit-level
Hausdorff consistency without further model-specific stochastic arguments.
\end{remark}

\bibliographystyle{plainnat}
\bibliography{references}

@article{LeoneNelsonNottingham1961,
  author  = {Leone, F. C. and Nelson, L. S. and Nottingham, R. B.},
  title   = {The Folded Normal Distribution},
  journal = {Technometrics},
  year    = {1961},
  volume  = {3},
  number  = {4},
  pages   = {543--550},
  doi     = {10.1080/00401706.1961.10489974}
}

@book{JohnsonKotzBalakrishnan1994,
  author    = {Johnson, N. L. and Kotz, S. and Balakrishnan, N.},
  title     = {Continuous Univariate Distributions, {V}ol.\ 1},
  edition   = {2},
  publisher = {Wiley},
  year      = {1994}
}

@article{NadarajahKotz2006,
  author  = {Nadarajah, S. and Kotz, S.},
  title   = {The Folded Normal Distribution},
  journal = {Statistical Papers},
  year    = {2006},
  volume  = {47},
  pages   = {545--558},
  
}

@book{vanderVaart1998,
  author    = {van der Vaart, A. W.},
  title     = {Asymptotic Statistics},
  publisher = {Cambridge University Press},
  year      = {1998}
}

@article{Teicher1961,
  author  = {Teicher, H.},
  title   = {Identifiability of Finite Mixtures},
  journal = {The Annals of Mathematical Statistics},
  year    = {1961},
  volume  = {32},
  number  = {1},
  pages   = {244--248},
  doi     = {10.1214/aoms/1177705155}
}

@article{Teicher1963,
  author  = {Teicher, H.},
  title   = {Identifiability of Mixtures},
  journal = {The Annals of Mathematical Statistics},
  year    = {1963},
  volume  = {34},
  number  = {4},
  pages   = {1265--1269},
  doi     = {10.1214/aoms/1177703862}
}

@article{YakowitzSpragins1968,
  author  = {Yakowitz, S. J. and Spragins, J. D.},
  title   = {On the Identifiability of Finite Mixtures},
  journal = {The Annals of Mathematical Statistics},
  year    = {1968},
  volume  = {39},
  number  = {1},
  pages   = {209--214},
  doi     = {10.1214/aoms/1177698520}
}

@article{AllmanMatiasRhodes2009,
  author  = {Allman, E. S. and Matias, C. and Rhodes, J. A.},
  title   = {Identifiability of Parameters in Latent Structure Models with Many Observed Variables},
  journal = {The Annals of Statistics},
  year    = {2009},
  volume  = {37},
  number  = {6A},
  pages   = {3099--3132},
  doi     = {10.1214/09-AOS689}
}

@article{HeinrichKahn2018,
  author  = {Heinrich, P. and Kahn, J.},
  title   = {Strong Identifiability and Optimal Rates for Finite Mixtures},
  journal = {The Annals of Statistics},
  year    = {2018},
  volume  = {46},
  number  = {6A},
  pages   = {2844--2870},
  doi     = {10.1214/17-AOS1641}
}

@article{Day1969,
  author  = {Day, N. E.},
  title   = {Estimating the Components of a Mixture of Normal Distributions},
  journal = {Biometrika},
  year    = {1969},
  volume  = {56},
  number  = {3},
  pages   = {463--474},
  doi     = {10.1093/biomet/56.3.463}
}

@article{Hathaway1985,
  author  = {Hathaway, R. J.},
  title   = {A Constrained Formulation of Maximum-Likelihood Estimation for Normal Mixture Distributions},
  journal = {The Annals of Statistics},
  year    = {1985},
  volume  = {13},
  number  = {2},
  pages   = {795--800},
  doi     = {10.1214/aos/1176349557}
}

@article{IngrassiaRocci2007,
  author  = {Ingrassia, S. and Rocci, R.},
  title   = {Constrained {EM} Algorithms for Mixture Models},
  journal = {Computational Statistics \& Data Analysis},
  year    = {2007},
  volume  = {51},
  number  = {11},
  pages   = {5339--5351},
  doi     = {10.1016/j.csda.2006.10.011}
}

@article{Keribin2000,
  author  = {Keribin, C.},
  title   = {Consistent Estimation of the Order of Mixture Models},
  journal = {Sankhy\={a}: The Indian Journal of Statistics, Series A},
  year    = {2000},
  volume  = {62},
  pages   = {49--66},
  doi     = {10.1214/aos/1013203454}
}

@book{vanderVaartWellner1996,
  author    = {van der Vaart, A. W. and Wellner, J. A.},
  title     = {Weak Convergence and Empirical Processes},
  publisher = {Springer},
  year      = {1996}
}

@incollection{NeweyMcFadden1994,
  author    = {Newey, W. K. and McFadden, D.},
  title     = {Large Sample Estimation and Hypothesis Testing},
  booktitle = {Handbook of Econometrics},
  editor    = {Engle, R. F. and McFadden, D. L.},
  volume    = {4},
  publisher = {Elsevier},
  year      = {1994},
  pages     = {2111--2245},
  doi       = {10.1016/S1573-4412(05)80006-6}
}

@book{AliprantisBorder2006,
  author    = {Aliprantis, C. D. and Border, K. C.},
  title     = {Infinite Dimensional Analysis: A Hitchhiker's Guide},
  edition   = {3},
  publisher = {Springer},
  year      = {2006}
}

@book{McLachlanPeel2000,
  author    = {McLachlan, G. J. and Peel, D.},
  title     = {Finite Mixture Models},
  publisher = {Wiley},
  year      = {2000},
  doi       = {10.1002/0471721182}
}

@book{Lindsay1995,
  author    = {Lindsay, B. G.},
  title     = {Mixture Models: Theory, Geometry and Applications},
  publisher = {Institute of Mathematical Statistics},
  year      = {1995}
}

@article{DempsterLairdRubin1977,
  author  = {Dempster, A. P. and Laird, N. M. and Rubin, D. B.},
  title   = {Maximum Likelihood from Incomplete Data via the {EM} Algorithm},
  journal = {Journal of the Royal Statistical Society: Series B (Methodological)},
  year    = {1977},
  volume  = {39},
  number  = {1},
  pages   = {1--38},
  doi     = {10.1111/j.2517-6161.1977.tb01600.x}
}

@article{RednerWalker1984,
  author  = {Redner, R. A. and Walker, H. F.},
  title   = {Mixture Densities, Maximum Likelihood and the {EM} Algorithm},
  journal = {SIAM Review},
  year    = {1984},
  volume  = {26},
  number  = {2},
  pages   = {195--239},
  doi     = {10.1137/1026034}
}
\end{document}